\documentclass[11pt]{article}
\usepackage[utf8]{inputenc}
\usepackage{amsmath,amssymb,amscd,epsfig,bbm}
\usepackage{stmaryrd,mathabx}
\usepackage{comment}
\usepackage{color}
\usepackage[T1]{fontenc}

\usepackage{interval}
\usepackage{float}

\usepackage[textsize=scriptsize,textwidth=2cm]{todonotes}
\usepackage{enumitem}
\usepackage{varwidth}
\setlist{nolistsep}
\usepackage[colorlinks]{hyperref}

\usepackage{mathtools}
\makeatletter
\DeclareRobustCommand\widecheck[1]{{\mathpalette\@widecheck{#1}}}
\def\@widecheck#1#2{%
    \setbox\z@\hbox{\m@th$#1#2$}%
    \setbox\tw@\hbox{\m@th$#1%
       \widehat{%
          \vrule\@width\z@\@height\ht\z@
          \vrule\@height\z@\@width\wd\z@}$}%
    \dp\tw@-\ht\z@
    \@tempdima\ht\z@ \advance\@tempdima2\ht\tw@ \divide\@tempdima\thr@@
    \setbox\tw@\hbox{%
       \raise\@tempdima\hbox{\scalebox{1}[-1]{\lower\@tempdima\box
\tw@}}}%
    {\ooalign{\box\tw@ \cr \box\z@}}}
\makeatother

\usepackage{amsthm}

\newtheorem{defi}{Definition}%[section]
\newtheorem{prop}[defi]{Proposition}
\newtheorem{theo}[defi]{Theorem}
\newtheorem{theofr}[defi]{Théorème}
\newtheorem{conj}[defi]{Conjecture}
\newtheorem{lemm}[defi]{Lemma}
\newtheorem{lemmfr}[defi]{Lemme}
\newtheorem{coro}[defi]{Corollary}

\theoremstyle{definition}
\newtheorem{rema}[defi]{Remark}
\newtheorem{exem}[defi]{Example}
\newtheorem{exems}[defi]{Examples}

\newcommand{\bdefi}{\begin{defi}}
\newcommand{\edefi}{\end{defi}}
\newcommand{\bprop}{\begin{prop}}
\newcommand{\eprop}{\end{prop}}
\newcommand{\btheo}{\begin{theo}}
\newcommand{\etheo}{\end{theo}}
\newcommand{\btheofr}{\begin{theofr}}
\newcommand{\etheofr}{\end{theofr}}
\newcommand{\blemm}{\begin{lemm}}
\newcommand{\elemm}{\end{lemm}}
\newcommand{\blemmfr}{\begin{lemmfr}}
\newcommand{\elemmfr}{\end{lemmfr}}
\newcommand{\brema}{\begin{rema}}
\newcommand{\erema}{\end{rema}}
\newcommand{\bexer}{\begin{exem}}
\newcommand{\eexer}{\end{exem}}
\newcommand{\bexems}{\begin{exems}}
\newcommand{\eexems}{\end{exems}}
\newcommand{\bconj}{\begin{conj}}
\newcommand{\econj}{\end{conj}}
\newcommand{\bcoro}{\begin{coro}}
\newcommand{\ecoro}{\end{coro}}

\usepackage{mathrsfs}
\renewcommand\mathcal{\mathscr}

\newcommand{\A}{{\cal A}}

\newcommand{\D}{{\cal D}}
\newcommand{\E}{{\cal E}}

\newcommand{\maths}[1]{{\mathbb #1}}  

\newcommand{\CC}{\maths{C}}

\newcommand{\HH}{\maths{H}}

\newcommand{\NN}{\maths{N}}

\newcommand{\RR}{\maths{R}}

\newcommand{\ZZ}{\maths{Z}}

\newcommand{\ra}{\rightarrow}
\newcommand{\bs}{\backslash}

\newcommand{\wt}[1]{{\widetilde{#1}}}

\newcommand{\ga}{\gamma}
\newcommand{\Ga}{\Gamma}

\newcommand{\bigO}{\operatorname{O}}

\newcommand{\card}{{\operatorname{Card}}}

\newcommand{\diag}{{\operatorname{diag}}}

\newcommand{\Leb}{\operatorname{Leb}}

\newcommand{\mBM}{m_{\rm BM}}

\newcommand\Perp{\operatorname{Perp}}

\newcommand{\ssm}{\!\smallsetminus\!}

\newcommand{\wtmBM}{\wt m_{\rm BM}}

\newcommand{\hdr}{{\HH}^2_\RR}

\newcommand{\PSL}{\operatorname{PSL}}

\newcommand{\PSLZ}{\operatorname{PSL}_{2}(\ZZ)}

\newcommand{\flow}[1]{{{\tt g}^{#1}}}  
\newcommand\normalout{\partial^1_{+}}
\newcommand\normalin{\partial^1_{-}}
\newcommand\normalpm{\partial^1_{\pm}}

\newcounter{fig}

\usepackage{ulem}

\title{Equidistribution of common perpendiculars \\ in negative
  curvature}
\author{Jouni Parkkonen and Frédéric Paulin}
\date{\today}

\begin{document}
\bibliographystyle{../alphanum}
\maketitle

\begin{abstract} 
Let $A^-$ and $A^+$ be properly immersed closed locally convex subsets
of a Riemannian manifold $M$ with pinched negative sectional
curvature.  When the Bowen-Margulis measure on $T^1M$ is finite and
mixing for the geodesic flow, we prove that the Lebesgue measures
along the common perpendiculars of length at most $t$ from $A^-$ to
$A^+$, counted with multiplicities and lifted to $T^1M$,
equidistribute to the Bowen-Margulis measure as $t\ra+\infty$.
When $M$ is locally symmetric with finite volume and the geodesic flow
is exponentially mixing, we give an error term for the asymptotic.  When
$T^1M$ is endowed with a bounded Hölder-continuous potential, and when
the associated equilibrium state is finite and mixing for the geodesic
flow, we prove the equidistribution of these Lebesgue measures weighted
by the amplitudes of the potential to the equilibrium state.
%ce dernier si extension avec potentiel
\end{abstract}

\section{Introduction}
\label{sec:intro}

Let $M$ be a nonelementary complete connected Riemannian good orbifold
with pinched sectional curvature at most $-1$, and let $(\flow
t)_{t\in\RR}$ be its geodesic flow on its unit tangent bundle
$T^1M$. Let $A^-$ and $A^+$ be proper nonempty properly immersed
closed locally convex subsets of $M$. A {\it common perpendicular}
from $A^-$ to $A^+$ is a locally geodesic path in $M$ starting
perpendicularly from $A^-$ and arriving perpendicularly to $A^+$ (see
\cite[\S 2.2]{ParPau14ETDS} or \cite[\S 2.4]{BroParPau19} for
explanations when the boundary of $A^-$ or $A^+$ is not smooth).  For
all $t> 0$, we denote by $\Perp (A^-, A^+, t)$ the set of common
perpendiculars from $A^-$ to $A^+$ with length at most $t$, considered
with multiplicities.  We refer to \cite[\S 3.3]{ParPau14ETDS} or
\cite[\S 12.1]{BroParPau19} for the definition of the multiplicities,
which are equal to $1$ if $M$ is a manifold and if $A^-$ and $A^+$ are
embedded and disjoint.  For every $\alpha\in \Perp (A^-, A^+, t)$, if
$\ell(\alpha)$ is its length and if $v^-_\alpha$ and $v^+_\alpha$ are
its initial and terminal unit tangent vectors, we denote by
$\Leb_\alpha$ the pushforward measure on $T^1M$, by the map $t\mapsto
\flow t v^-_\alpha $, of the Lebesgue measure on $[0,\ell(\alpha)]$.

We denote the total mass of any measure $m$ by $\|m\|$. Referring to
\cite{Roblin03} (see also Section \ref{sec:geometry}) for definitions,
we denote by $\delta$ the critical exponent of the fundamental group
of $M$ and by $\mBM$ a Bowen-Margulis measure on $T^1M$. By
\cite{OtaPei04} and \cite{DilTho25}, if $\mBM$ if finite, then
$\delta$ is the topological entropy of the geodesic flow $(\flow
t)_{t\in\RR}$ and $\frac{\mBM}{\|\mBM\|}$ is its unique measure of
maximal entropy. We denote by $\sigma^\mp_{A^\pm}$ the skinning
measures of $A^\pm$ defined in \cite{ParPau14ETDS} (see also Section
\ref{sec:geometry}), generalising \cite{OhSha12,OhSha13} when $M$ has
constant curvature and $A^-,A^+$ are balls, horoballs or totally
geodesic submanifolds.

\btheo \label{theo:mainintro} Assume that the Bowen-Margulis measure
$\mBM$ is finite and mixing for the geodesic flow of $M$. Assume that
the skinning measures $\sigma^+_{A^-}$ and $\sigma^{-}_{A^+}$ are
finite and nonzero. For the narrow convergence of measures on $T^1M$,
we have
\begin{equation}\label{eq:mainintro}
\lim_{t\ra+\infty}\; \frac{\delta\;\|\mBM\|}
    {t\;e^{\delta\, t}\;\|\sigma^+_{A^-}\|\;\|\sigma^{-}_{A^+}\|}
\sum_{\alpha\in\Perp(A^-,\,A^+,\,t)} \; 
\Leb_\alpha =\frac{\mBM}{\|\mBM\|}\,.
\end{equation}
If furthermore $M$ is locally symmetric with finite volume, and if the
geodesic flow of $M$ is exponentially mixing, then there exists
$\ell\in\NN$ such that for every compact subset $K$ in $T^1M$ and every
$C^\ell$-smooth function $\psi:T^1M\ra\CC$ with support in $K$ and
$W^{\ell,2}$-Sobolev norm $\|\psi\|_{\ell}$, as $s\ra+\infty$, we
have
\[\frac{\delta\;\|\mBM\|}
{t\;e^{\delta\, t}\;\|\sigma^+_{A^-}\|\;\|\sigma^{-}_{A^+}\|}
\sum_{\alpha\in\Perp(A^-,\,A^+,\,t)} \; \Leb_\alpha(\psi)
=\frac{\mBM(\psi)}{\|\mBM\|} +\bigO_K(\frac{1}{t}\|\psi\|_{\ell})\,.
\]
\etheo

This result is an analog of the theorems of Bowen and Margulis on the
equidistribution towards the Bowen-Margulis measure of the Lebesgue
measures along the periodic orbits of the geodesic flows when $M$ is
compact (see \cite[Theo.~5.1.1]{Roblin03} for the extension to the
above generality, and \cite[\S 9.3]{PauPolSha15} for the extension to
equilibrium states).

As a very special case which is already new and striking, if we choose
$A^-=A^+=\{p\}$ for any $p\in M$, Theorem \ref{eq:mainintro} shows
that the Lebesgue mesures along geodesic loops based at $p$, lifted to
$T^1M$, equidistribute to the Bowen-Margulis measure as the upper
bound on their lengths increases to $+\infty$.

\noindent\begin{minipage}{8.7cm}
\begin{figure}[H]
\begin{center}
\includegraphics[width=2.6cm]{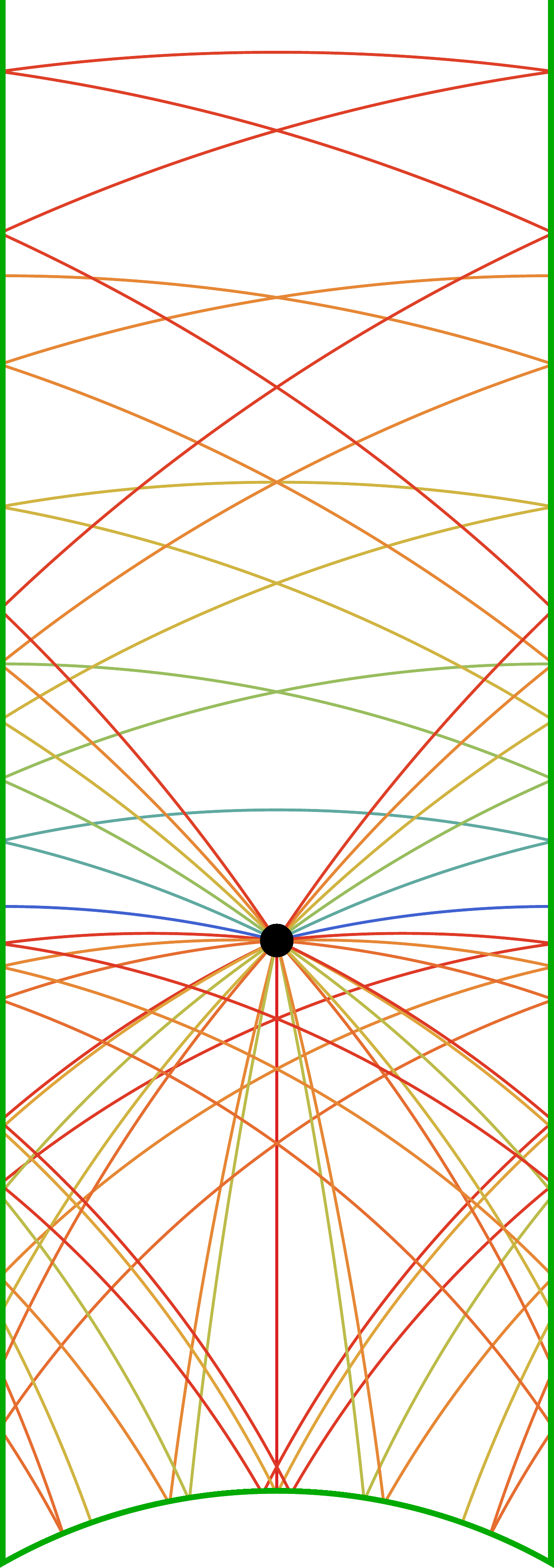}
\hspace{.1cm}
\includegraphics[width=2.6cm]{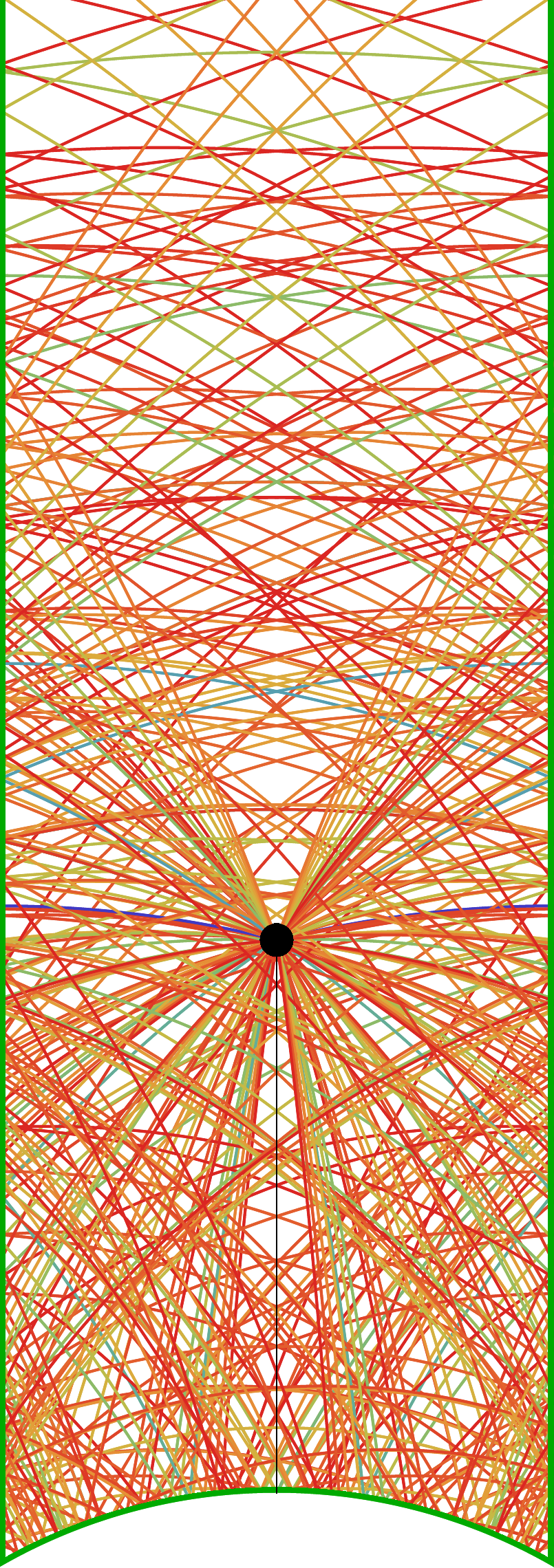}
\hspace{.1cm}
\includegraphics[width=2.6cm]{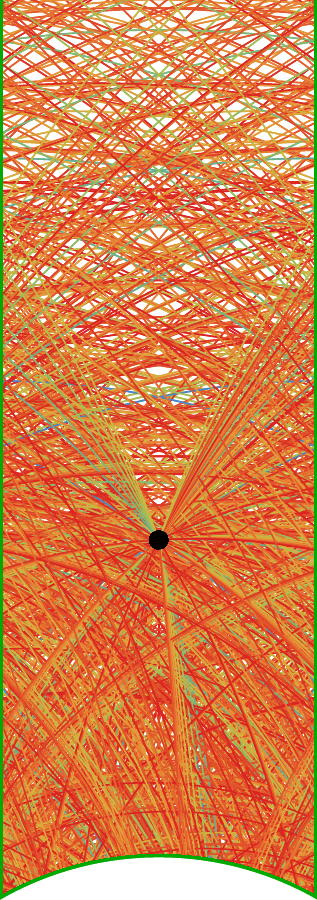}
\hspace{.1cm}
  \end{center}\vspace*{-0.7cm}
  \caption{Geodesic loops equidistribute.}
\end{figure}
\end{minipage}
\begin{minipage}{6.2cm}
~~ The figures on the left show geodesic loops based at the image of the point
$2i\in\hdr$ in the modular orbifold $M=\PSL_2(\ZZ)\bs\hdr$, with
$\hdr$ the upper halfspace model of the real hyperbolic plane.  These
geodesic loops, that are in general not closed geodesics, are shown
lifted to the usual fundamental domain of $\PSLZ$ with its boundary
identifications.  The first figure displays the geodesic loops of
length at most $2.5$, and the second and third figures display the geodesic
loops of length at most $4$ and $5$ restricted to the same part of $M$.
See also \cite{LiSta23} and \cite{PolWar} for related recent work.
\end{minipage}

\bigskip
Theorem  \ref{theo:mainintro} builds on the joint equidistribution theorem
\cite[Theo.~2]{ParPau17ETDS} of the authors, stated as Theorem
\ref{theo:17ETDS2} in Section \ref{sec:proof}, for the pairs
$(v^-_\alpha, v^+_\alpha)$ of initial and terminal unit tangent
vectors of the elements $\alpha\in\Perp (A^-, A^+, t)$ in the
outer/inner normal bundles of $A^\pm$. But Theorems
\ref{theo:mainintro} and \ref{theo:17ETDS2} are very different in
nature, the second one being an equidistribution on $T^1M\times T^1M$
towards the product of two measures whose supports have zero measure
for the Bowen-Margulis measure on $T^1M$. We refer to
\cite{ParPau14ETDS,ParPau17ETDS,ParPau17MA,BroParPau19,ParPau22MPCPS,
  ParPau25}, and in particular to the surveys \cite{ParPau16LMS,
  ParPau17CIRM}, for other motivations, arithmetic applications and
references for the study of common perpendiculars in negative
curvature.
%\end{minipage}

We have the following version of Theorem \ref{theo:mainintro} for
equilibrium states. We refer to \cite{PauPolSha15,DilTho25} and
Subsection \ref{subsec:thermoform} for the various definitions of the
objects in the following statement.

\bigskip

\btheo \label{theo:intropot} Let $F:T^1M\ra \CC$ be a bounded
Hölder-continuous function with positive topological pressure
$\delta_F$ and amplitudes $\int_\alpha F$ along every $\alpha\in
\Perp(A^-,\,A^+, \,t)$.  Assume that the Gibbs measure $m_F$ of $F$ is
finite and mixing for the geodesic flow of $M$. Assume that the
skinning measures $\sigma^+_{A^-}$ and $\sigma^{-}_{A^+}$ associated
with $F$ are finite and nonzero. For the narrow convergence of
measures on $T^1M$, we have
\[
\lim_{t\ra+\infty}\; \frac{\delta_F\;\|m_F\|}
    {t\;e^{\delta_F\, t}\;\|\sigma^+_{A^-}\|\;\|\sigma^{-}_{A^+}\|}
\sum_{\alpha\in\Perp(A^-,\,A^+,\,t)} \; 
e^{\int_\alpha F}\Leb_\alpha =\frac{m_F}{\|m_F\|}\,.
\]
\etheo

The first claim of Theorem \ref{theo:mainintro} is the particular case
when $F=0$ of Theorem \ref{theo:intropot}. But for the convenience of
the reader, we will give a full proof of Theorem \ref{theo:mainintro}
in Section \ref{sec:proof}, and only indicate in Section
\ref{sec:potential} the quite technical changes brought to that proof
by the potential $F$.

Theorem \ref{theo:mainintro} will be an important tool in two
forthcoming papers. In \cite{ParPauSay25} we study the
equidistribution of divergent geodesics in $M$ (when $M$ has finite
volume), and in \cite{ErlParPau24} we study the equidistribution of
reciprocal closed geodesics in $M$.

\section{The geometric and measure theoretic background}
\label{sec:geometry}

We refer to \cite{BriHae99} for background on the geometric content of
this section. Let $\wt M$ be a complete simply connected Riemannian
manifold with (dimension at least $2$ and) pinched negative sectional
curvature $-b^2\le K\le -1$, and let $x_*\in\wt M$ be a fixed
basepoint.

Let $\Ga$ be a nonelementary (not virtually nilpotent) discrete group
of isometries of $\wt M$, let $M$ be the quotient (hence good)
Riemannian orbifold $\Ga\bs\wt M$, and let $T^1M$ be the quotient
orbifold $\Ga\bs T_1\wt M$.  We denote by $\partial_{\infty}\wt M$ the
boundary at infinity of $\wt M$, by
\[
\delta= \lim_{t\ra+\infty}
\frac{1}{t} \,\card\,\{\ga\in\Ga\,:\, d(x_*,\ga x_*)\leq t\}
\]
the critical exponent of $\Ga$, and by $\Lambda\Ga$ the limit set of
$\Ga$.  We denote by $\pi$ the footpoint projections $T^1\wt M\ra\wt
M$ and $T^1M \ra M$.  Let $(\flow t)_{t\in\RR}$ be the geodesic flow
on the unit tangent bundle $T^1\wt M$ of $\wt M$, as well as its
quotient flow on $T^1M$.

For every $v\in T^1\wt M$, let $v_-\in\partial_\infty\wt M$ and
$v_+\in\partial_\infty\wt M$, respectively, be the endpoints at
$-\infty$ and $+\infty$ of the geodesic line defined by $v$. We denote
by $p_\pm:T^1\wt M\ra \partial_\infty\wt M$ the {\it endpoint maps}
$v\mapsto v_\pm$.  Let us denote by $\diag=\{(x,x):x\in X\}$ the
diagonal in any Cartesian square $X\times X$ of a set $X$. {\it Hopf's
  parametrisation} with respect to the point $x_*$ of $T^1\wt M$
is the homeomorphism which identifies $T^1\wt M$ with
$(\partial_\infty \wt M\times\partial_\infty\wt M\ssm\diag)\times\RR$
by the map $v\mapsto(v_-,v_+,s)$, where $s$ is the signed distance to
$\pi(v)$ of the closest point to $x_*$ on the geodesic line defined by
$v$.

For every $\xi\in\partial_{\infty}\wt M$, let $\rho_{\xi}:
[0,+\infty[\; \ra \wt M$ be the geodesic ray with origin $x_{*}$ and
point at infinity $\xi$. The {\it Busemann cocycle} of $\wt M$ is the
map $\beta: \wt M\times\wt M\times\partial_{\infty} \wt M\to\RR$
defined by
\begin{equation}\label{eq:defiBuse}
  (x,y,\xi)\mapsto \beta_{\xi}(x,y)= \lim_{t\to+\infty}
  d(\rho_{\xi}(t),x)-d(\rho_{\xi}(t),y)\,.
\end{equation}
The {\it visual distance} $d_{x_*}$ on $\partial_\infty \wt M$ seen
from $x_*$ is defined by $d_{x_*}(\xi,\eta)= e^{-\frac{1}{2}
  (\beta_\xi(x_*,\,y)+\beta_\eta(x_*,\,y))}$ where $y$ is the closest
point to $x_*$ on the geodesic line between two distinct points at
infinity $\xi$ and $\eta$.

\medskip
We refer to \cite{Roblin03} for more background and for the basic
properties of the following measures. A family $(\mu_x)_{x\in\wt M}$
of finite measures on $\partial_\infty\wt M$, whose support is the
limit set $\Lambda\Ga$ of $\Ga$, is a {\it Patterson-Sullivan density}
for $\Ga$ if
\[
\ga_*\mu_x=\mu_{\ga x}
\]
for all $\ga\in\Ga$ and $x\in \wt M$, and if the following
Radon-Nikodym derivatives exist for all $x,y\in\wt M$ and satisfy for
(almost) every $\xi\in\partial_\infty\wt M$
\[
\frac{d\mu_x}{d\mu_y}(\xi)=e^{-\delta\,\beta_\xi(x,\,y)}\,.
\]

We fix such a family $(\mu_x)_{x\in\wt M}$.  The {\it Bowen-Margulis
  measure} on $T^1\wt M$ (associated with this Patterson-Sullivan
density) is the measure $\wtmBM$ on $T^1\wt M$ given by the density
\begin{equation}\label{eq:defimBM}
d\wtmBM(v)=
e^{-\delta(\beta_{v_-}(\pi(v),\,x_*)\,+\,\beta_{v_+}(\pi(v),\,x_*))}\;
d\mu_{x_*}(v_-)\,d\mu_{x_*}(v_+)\,dt
\end{equation} 
in Hopf's parametrisation of $T^1\wt M$ with respect to $x_*$. The
Bowen-Margulis measure $\wtmBM$ is independent of $x_*$, and it is
invariant under the actions of the group $\Ga$ and of the geodesic
flow. Thus, it defines a measure $\mBM$ on $T^1M$ which is invariant
under the quotient geodesic flow, called the {\it Bowen-Margulis
  measure} on $T^1M$.  If $\mBM$ is finite, then the
Patterson-Sullivan densities are unique up to a multiplicative
constant; hence the Bowen-Margulis measure is uniquely defined, up to
a multiplicative constant.

Babillot \cite[Theo.~1]{Babillot02} showed that if $\mBM$ is finite,
then it is mixing for the geodesic flow of $M$ if the length spectrum
of $M$ is not contained in a discrete subgroup of $\RR$, as in
particular when $M$ is locally symmetric.  For every Riemannian
orbifold $Y$, we denote by $\|\cdot\|_{\ell}$ the
$W^{\ell,2}$-Sobolev norm on the vector space $C^{\ell}_c(Y)$ of
$C^\ell$-smooth functions with compact support on $Y$.  We refer for
instance to \cite[\S 9.1]{BroParPau19} for the definition of the
exponentially mixing property (for the Sobolev regularity) of the
geodesic flow of $M$.  Note that when $M$ is locally symmetric with
finite volume, by the work of Li-Pan \cite{LiPan22} when $M$ is real
hyperbolic and by the Margulis arithmeticity result with the works of
Kleinbock-Margulis and Clozel, see for instance \cite[page
  182]{BroParPau19}, when $M$ is quaternionic or octonionic
hyperbolic, the only case when the geodesic flow of $M$ is not yet
known to be exponentially mixing is when $M$ is complex hyperbolic.

\medskip
Let $D$ be a nonempty proper closed convex subset of $\wt M$. We refer
to \cite[\S 2.2]{ParPau14ETDS} or \cite[\S 2.4]{BroParPau19} for the
definition of the outer and inner normal bundles $\normalout D$ and
$\normalin D$ of $D$. We refer to \cite{ParPau14ETDS} for more
background and for the basic properties of the following measures. The
{\it (outer) skinning measure} on $\normalout D$ (associated with the
Patterson-Sullivan density $(\mu_x)_{x\in\wt M}$) is the measure
$\wt\sigma^+_D$ on $\normalout D$ defined, using the positive endpoint
homeomorphism $p_+:v\mapsto v_+$ from $\normalout D$ to
$\partial_\infty\wt M-\partial_\infty D$, by
\[
d\,\wt\sigma^+_D(v) = e^{-\delta\,
  \beta_{v_+}(\pi(v),\,x_*)}\;d\mu_{x_*}(v_{+}) \,.
\]
The {\it (inner) skinning measure} $d\wt\sigma^-_D(v)= e^{-\delta\,
  \beta_{v_-}(\pi(v),\,x_*)}\; d\mu_{x_*}(v_{-})$ is the similarly
defined measure on $\normalin D$.  When $D$ is a singleton, we
immediately have
\begin{equation}\label{eq:skinningsingle}
  \forall\;v\in T^1_{x_*}\wt M,\qquad d\,\wt \sigma_{\{x_*\}}^\pm(v)=
  d\mu_{\{x_*\}}(v_\pm)\,.
\end{equation}
If the family $(\ga D)_{\ga \in\Ga/\Ga_{D}}$ is locally finite in $\wt
M$, and if $A$ is the image of $D$ in $M$, we say that $A$ is a {\it
  proper nonempty properly immersed closed locally convex subset} of
$M$. We denote by $\sigma^\pm_A$ the locally finite Borel measure on
$T^1 M$ induced by the $\Ga$-invariant locally finite Borel measure
$\sum_{\ga \in\Ga/\Ga_{D}} \ga_*\,\wt\sigma^\pm_D$ on $T^1 \wt M$. The
support of $\sigma^\pm_A$ is the image $\normalpm A$ of $\normalpm D$
by the map $T^1 \wt M\ra T^1M$.

\medskip We conclude this section by the following new topological
construction.  In the same way the geometric compactification
$\overline{\wt M}=\wt M\cup \partial_\infty \wt M$ of $\wt M$
compactifies $\wt M$ by gluing its boundary at infinity
$\partial_\infty \wt M$ (see \cite{BriHae99}), there exists a {\it
  geometric compactification} $\overline{T^1\wt M}= T^1\wt M\cup
\partial_\infty \wt M$ of $T^1\wt M$ by gluing $\partial_\infty \wt
M$, in which $T^1\wt M$ is open and dense, constructed as follows. The
footpoint projection $\pi: T^1\wt M\ra \wt M$ uniquely extends by the
identity map on $\partial_\infty \wt M$ to a set-theoretic map
$\overline{T^1\wt M}\ra \overline{\wt M}$ again denoted by $\pi$. A
basis of open subsets for the topology of $\overline{T^1\wt M}$
consists of either the open subsets of $T^1\wt M$ or of the preimages
by $\pi$ of the open subsets of $\overline{\wt M}$.

It is well known that the uniform structure\footnote{See \cite[Chap.~2]
{Bourbaki71a} for background on uniform spaces.} on $\partial_\infty
\wt M$ defined by the visual distance $d_{x_*}$ extends to a uniform
structure on $\overline{\wt M}$ (hence by pullback by $\pi$ on
$\overline{T^1\wt M}$), as follows. Recall (see \cite{Bourdon95}) that
the visual distance between two distinct points at infinity
$\xi,\eta\in \partial_\infty M$ is comparable to the exponential of
the opposite of the distance from $x_*$ to the geodesic line between
$\xi$ and $\eta\,$: There exists a universal constant $c\geq 0$ such that
$e^{-d(x_*,\,\interval[open]{\xi}{\eta}\,)}\leq d_{x_*}(\xi, \eta)
\leq e^{-d(x_*,\,\interval[open]{\xi}{\eta}\,) +c}$.  A fundamental
system of entourages $(\E_\epsilon)_{\epsilon>0}$ for the uniform
structure on $\overline{\wt M}$ consists of the sets $\E_\epsilon$ of
pairs $(x,y)$ in $\overline{\wt M} \times \overline{\wt M}$ such that
the distance between $x_*$ and the geodesic segment, ray or line
between $x$ and $y$ is at least $-\ln \epsilon$.

\section{Equidistribution of Lebesgue measures along
  common perpendiculars}
\label{sec:proof}

This whole Section is devoted to the proof of Theorem
\ref{theo:mainintro}.  We fix $D^-$ and $D^+$ two nonempty proper
closed convex subsets of $\wt M$, such that the families $(\alpha^\pm
D^\pm)_{\alpha^\pm \in \Ga/\Ga_{D^\pm}}$ are locally finite in $\wt
M$, and we denote by $A^-$ and $A^+$ respectively their images in $M$.
Using the notation of the introduction $v^\pm_\alpha$ for
$\alpha\in\Perp(A^-,\,A^+, \,t)$, one of the tools of the proof of
Theorem \ref{theo:mainintro} is the following result. It is proved in
\cite[Coro.~12]{ParPau17ETDS} for its first claim and \cite[Theo.~15
  (2)]{ParPau17ETDS} for its second claim (with an extra smoothing
argument of the boundary of $D^\pm$).

\btheo\label{theo:17ETDS2} Assume that the Bowen-Margulis measure
$\mBM$ is finite and mixing for the geodesic flow of $M$. Assume that
the skinning measures $\sigma^+_{A^-}$ and $\sigma^{-}_{A^+}$ are
finite and nonzero. For the narrow convergence of measures on
$T^1M\times T^1M$, we have
\[
\;\;\;\lim_{t\ra+\infty}\; \delta\;\|\mBM\|\;e^{-\delta\, t}
\sum_{\alpha\in\Perp(A^-,\,A^+,\,t)} \; \Delta_{v^-_\alpha}
\otimes\Delta_{v^+_\alpha}\;=\; \sigma^+_{A^-}\otimes
\sigma^-_{A^+}\,.
\]
If furthermore $M$ is locally symmetric with finite volume, and if the
geodesic flow of $M$ is exponentially mixing, then there exists
$\ell\in\NN$ and $\kappa>0$ such that for all $\phi^\pm\in C^\ell_c
(T^1M)$, as $s\ra+\infty$, we have
\[
\frac{\delta\;\|\mBM\|}{e^{\delta\, t}}
\sum_{\alpha\in\Perp(A^-,\,A^+,\,t)} \;
\phi^-(v^-_\alpha)\,\phi^+(v^+_\alpha)
=\sigma^+_{A^-}(\phi^-)\;\sigma^-_{A^+}(\phi^+)+ \bigO(e^{-\kappa
  t}\,\|\phi^-\|_{\ell}\,\|\phi^+\|_{\ell})\,.\quad\Box
\]
\etheo

We start the proof of Theorem \ref{theo:mainintro} by giving the
notation that we will use.  Let $D'$ and $D''$ be nonempty closed
convex subsets of $\wt M$ such that $d(D',D'')>0$. Let $[D',D'']$ be
the common perpendicular arc between $D'$ and $D''$ (oriented from
$D'$ towards $D''$). Let $u_{D',D''}$ (respectively $v_{D',D''}$) be
the initial (respectively terminal) unit tangent vector to $[D',D'']$,
so that $v_{D',D''}= \flow{d(D',D'')} u_{D',D''}$. Let
$\Leb_{[D',D'']}$ be the measure on $T^1\wt M$ which is the
pushforward measure by the continuous map $t\mapsto
\flow{t}u_{D',D''}$ of the Lebesgue measure on $[0,d(D',D'')]$.  Its
support is
\[
T^1[D',D'']=\big\{\flow{t}u_{D',D''}: t\in [0,d(D',D'')]\big\}\,.
\]
For every isometry $\ga$ of $\wt M$, we naturally have $\ga
u_{D',D''}= u_{\ga D',\ga D''}$, $\ga v_{D',D''} =v_{\ga D',\ga D''}$
and $\ga_* \Leb_{[D',D'']}= \Leb_{[\ga D',\ga D'']}$.

For every $t>0$, let
\[
\nu_{1,t}=\sum_{\alpha^-\in\Ga/\Ga_{D^-},\,
  \alpha^+\in\Ga/\Ga_{D^+}\,:\; 0<d( \alpha^-D^-,\,\alpha^+D^+)\leq t}
\Leb_{[\alpha^-D^-,\alpha^+D^+]}\,,
\]
which is a $\Ga$-invariant Borel measure on $T^1\wt M$. We claim that
this measure $\nu_{1,t}$ is locally finite. Indeed, let $K$ be a
nonempty compact subset of $T^1\wt M$, let $t_K$ be the diameter of
the subset $\pi(K)$ of $\wt M$, and let $x_K\in \pi(K)$. For all
$\alpha^\pm\in \Ga/\Ga_{D^\pm}$ such that $K$ meets the support of
$\Leb_{[\alpha^-D^-,\alpha^+D^+]}$, we have $d(x_K,\alpha^\pm D^\pm)
\leq t+t_K$. Since the families $(\alpha^\pm D^\pm)_{\alpha^\pm
  \in\Ga/\Ga_{D^\pm}}$ are locally finite, the sets
\[
E^\pm=\{\alpha^\pm\in\Ga/\Ga_{D^\pm}:d(x_K,\alpha^\pm D^\pm) \leq
t+t_K\}
\]
are finite. Hence ${\nu_{1,t}}_{\,\mid K}=\sum_{\alpha^\pm\in
  E^\pm\,:\; 0<d( \alpha^-D^-,\,\alpha^+D^+)\leq t}
{\Leb_{[\alpha^-D^-,\,\alpha^+D^+]}} _{\,\mid K}$ is a finite sum of
finite measures on $K$.

The locally finite $\Ga$-invariant Borel measure $\nu_{1,t}$ induces
by the orbifold covering map $T^1\wt M\ra T^1 M=\Ga\bs T^1 \wt M$ a
locally finite Borel measure on $T^1 M$, see \cite[\S 2.6]
{PauPolSha15}. This measure is exactly
$\sum_{\alpha\in\Perp(A^-,A^+,t)} \Leb_\alpha$, which is, up to
normalization, the measure in the left hand-side of Equation
\eqref{eq:mainintro}. To see this, note that the Lebesgue measure
$\Leb_\alpha$ along every $\alpha\in\Perp(A^-,A^+,t)$ lifts to the
Lebesgue measure $\Leb_{[\alpha^-D^-,\alpha^+D^+]}$ along a common
perpendicular between the images (at distance at most $t$) of $D^-$
and $D^+$ by two elements $\alpha^-$ and $\alpha^+$ of $\Ga$ well
defined modulo the stabilizers of $D^-$ and $D^+$, and conversely. The
multiplicities are defined in \cite[\S 3.3]{ParPau17ETDS}, see also
\cite[\S 12.1]{BroParPau19}, precisely in order to deal with the
orbifold covering and with the multiple intersections of the locally
finite families $(\alpha^\pm D^\pm)_{\alpha^\pm\in\Ga/\Ga_{\wt D^\pm}}$.

Let us fix $\epsilon \in \interval[open]{0}{\frac{1}{2}}$ and a point
$x_0\in\wt M$.  At the very end of the proof, $\epsilon$ will tend to
$0$ and $x_0$ will vary in an $\epsilon$-net of $\wt M$. We will
denote by $\epsilon_1,\epsilon_2,\epsilon_3$ positive functions of
$\epsilon$ (depending only on $(\wt M,\Ga)$) that converge to $0$ as
$\epsilon$ tends to $0$. Let $B=\pi^{-1}(B(x_0, \frac{\epsilon}{2}))$
be the subset of $T^1\wt M$ of elements whose footpoints are at
distance less than $\frac{\epsilon}{2}$ from $x_0$.  Let $\psi:T^1\wt
M\ra\RR$ be a nonzero nonnegative continuous function with compact
support contained in $B$. Using Hopf's parametrisation with respect to
the given point $x_0$, that identifies $T^1\wt M$ with
$(\partial_\infty \wt M\times \partial_\infty\wt M\ssm\diag)\times
\RR$, we also assume that $\psi$ is a product of three nonzero
nonnegative continuous functions with compact support
$\psi^-:\partial_\infty\wt M\ra\RR$, $\psi^+:\partial_\infty\wt M \ra
\RR$ and $\psi^0:\partial_\infty\wt M\ra\RR$, on each factor of the
product $\partial_\infty\wt M\times\partial_\infty\wt M\times \RR$, so
that
\[
\psi:(\xi,\eta,s)\mapsto \psi^-(\xi)\;\psi^+(\eta)\;\psi^0(s)\,.
\]

The following picture will be useful throughout the proof.

\begin{figure}[H]
  \begin{center}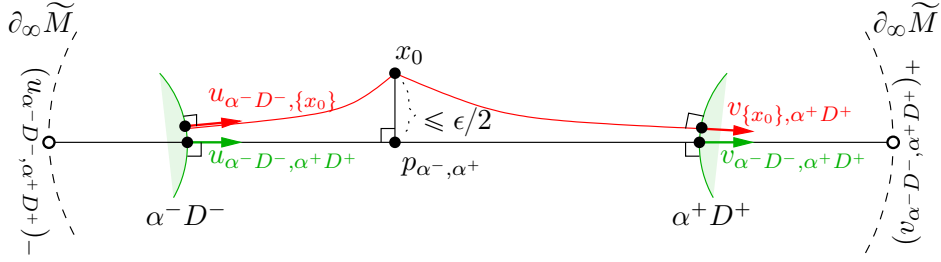\end{center}
  \vspace*{-0.7cm}
  \caption{\label{fig:perpcomm} Splitting common perpendiculars.}
\end{figure}

For every $t>0$, using Hopf's parametrisation with respect to $x_0$,
we have
\[
\nu_{1,t}=\sum_{\stackrel{\alpha^-\in\Ga/\Ga_{D^-},\,
  \alpha^+\in\Ga/\Ga_{D^+}}{0<d( \alpha^-D^-,\,\alpha^+D^+)\leq t}}
\big(\Delta_{(u_{\alpha^-D^-,\alpha^+D^+})_-}\otimes
\Delta_{(v_{\alpha^-D^-,\alpha^+D^+})_+}\otimes
ds\big)_{\,\mid\, T^1[\alpha^-D^-,\alpha^+D^+]}\,.
\]
For every $t>0$, let us consider the $\Ga$-invariant Borel measure
$\nu_{2,t}$ on $(T^1\wt M\times T^1\wt M\ssm\diag)\times \RR$ defined
by
\[
\nu_{2,t}=\sum_{\stackrel{\alpha^-\in\Ga/\Ga_{D^-},\,
  \alpha^+\in\Ga/\Ga_{D^+}}{0<d( \alpha^-D^-,\,\alpha^+D^+)\leq t}}
\Delta_{u_{\alpha^-D^-,\alpha^+D^+}}\otimes
\Delta_{v_{\alpha^-D^-,\alpha^+D^+}}\otimes ds\,.
\]
The measure $\nu_{2,t}$ is locally finite since the set of initial
unit tangent vectors $u_{\alpha^-D^-,\alpha^+D^+}$ for $\alpha^\pm\in
\Ga/ \Ga_{D^\pm}$ such that $0<d( \alpha^-D^-,\,\alpha^+D^+)\leq t$
and the corresponding set of terminal unit tangent vectors
$v_{\alpha^-D^-, \alpha^+D^+}$ are locally finite in $T^1\wt M$ , and
the support of $\nu_{2,t}$ is the product of these two locally finite
sets and $\RR$.  Both measures $\nu_{1,t}$ and $\nu_{2,t}$ can be seen
as measures on $(\overline{T^1\wt M}\times \overline{T^1\wt
  M}\ssm\diag)\times \RR$ using the geometric compactification of
$T^1\wt M$ defined at the end of Section \ref{sec:geometry}.

We endow $T^1\wt M=(\partial_\infty\wt M\times \partial_\infty\wt
M\ssm\diag) \times \RR$ with a product distance of the visual distance
$d_{x_0}$ seen from $x_0$ on each factor $\partial_\infty\wt M$ and
the Euclidean distance on the factor $\RR$.  We extend the function
$\psi=\psi^-\times\psi^+\times\psi^0$, for now defined on
$(\partial_\infty\wt M\times \partial_\infty \wt M\ssm\diag) \times
\RR$, to a continuous function on $(\overline{T^1\wt M}\times
\overline{T^1\wt M} \ssm\diag)\times \RR$, which is a product of three
continuous func\-tions again denoted by $\psi^-$, $\psi^+$, $\psi^0$
on each factor of the product $\overline{T^1\wt M}\times
\overline{T^1\wt M} \times \RR$.  We may assume that the supports of
the extended functions $\psi^\pm$ are contained in the complementary
subset of $\pi^{-1}(B(x_0,-\ln \epsilon))$ without changing their
modulus of continuity (up to a constant) for the uniform structure.

Let $\alpha^-\in\Ga/\Ga_{D^-}$ and $\alpha^+\in \Ga/\Ga_{D^+}$ be
elements that give a nontrivial contribution to the sum $\nu_{2,t}
(\psi)$. We then have $x_0\notin \alpha^-D^- \cup\alpha^+D^+$, since
$d(x_0, \pi(u_{\alpha^-D^-,\alpha^+D^+})\geq -\ln\epsilon$, hence
$d(x_0, \alpha^-D^-)\geq -\ln\epsilon-\frac{\epsilon}{2}\geq
-\ln\frac{1}{2} -\frac{1}{4}>0$, and similarly $d(x_0,\alpha^+D^+)
>0$. Furthermore, the vectors $u_{\alpha^-D^-,\alpha^+D^+}$ and
$v_{\alpha^-D^-,\alpha^+D^+}$ are close in $\overline{T^1\wt M}$ to
the points at infinity $(u_{\alpha^-D^-,\alpha^+D^+})_-$ and
$(v_{\alpha^-D^-,\alpha^+D^+})_+$ respectively, uniformly in such
$\alpha^-$ and $\alpha^+$.  By the uniform continuity of $\psi$, for
$t$ large enough, we have
\begin{equation}\label{eq:nu1}
e^{-\epsilon_1}\nu_{1,t}(\psi) \leq \nu_{2,t}(\psi)\leq
e^{\epsilon_1}\nu_{1,t}(\psi)\,.
\end{equation}
The vectors $u_{\alpha^-D^-,\alpha^+D^+}$ and
$v_{\alpha^-D^-,\alpha^+D^+}$ are uniformly close in $T^1\wt M$ to the
vectors $u_{\alpha^-D^-,\{x_0\}}$ and $v_{\{x_0\},\alpha^+D^+}$
respectively, which themselves are uniformly close in
$\overline{T^1\wt M}$ to the points at infinity
$(u_{\alpha^-D^-,\{x_0\}})_-$ and $(v_{\{x_0\},\alpha^+D^+})_+$. Thus
similarly, if
\[
\nu_{3,t}=\sum_{\stackrel{\alpha^-\in\Ga/\Ga_{D^-},\,
    \alpha^+\in\Ga/\Ga_{D^+}}{x_0\,\notin\, \alpha^-D^-
\cup\alpha^+D^+,\; 0<d( \alpha^-D^-,\,\alpha^+D^+)\leq t}}
\Delta_{(u_{\alpha^-D^-,\{x_0\}})_-}\otimes
\Delta_{(v_{\{x_0\},\alpha^+D^+})_+}\otimes ds\,,
\]
then $\nu_{3,t}$ is a $\Ga$-invariant Borel measure on $T^1\wt
M=(\partial_\infty\wt M\times \partial_\infty\wt M\ssm\diag) \times
\RR$ such that
\begin{equation}\label{eq:nu2}
e^{-\epsilon_2}\nu_{2,t}(\psi) \leq \nu_{3,t}(\psi)\leq
e^{\epsilon_2}\nu_{2,t}(\psi)\,.
\end{equation}

For every $t> 1$, let $N=\lceil\frac{t}{\epsilon}\rceil$, so that
$(N-1)\epsilon<t\leq N\epsilon$. For every $k\in\NN$, let
\[
\A_{k}=\big\{(\alpha^-,\alpha^+)\in\Ga/\Ga_{D^-}\times\Ga/\Ga_{D^+}:
x_0\notin \alpha^\pm D^\pm,\;
\max\{0,k\epsilon-1\}<d(\alpha^-D^-,\alpha^+D^+) \leq k\epsilon\big\}
\]
and
\begin{equation}\label{eq:defiSigmak}
Z_k=
\sum_{(\alpha^-,\alpha^+)\in\A_{k}}\psi^-\big((u_{\alpha^-D^-,\{x_0\}})_-\big)\;
  \psi^+\big((v_{\{x_0\},\alpha^+D^+})_+\big)\int\psi^0\, ds\,.
\end{equation}
We then have
\begin{equation}\label{eq:relatnu3tSigmak}
  \sum_{k=1}^{N-1}\;Z_k\leq \nu_{3,t}(\psi)\leq
  \sum_{k=1}^N \;Z_k\,.
\end{equation}
For all $i,j\in\ZZ$, let us define
\[
\A^-_i=\big\{\alpha^- \in \Ga/\Ga_{D^-}:
\max\{0,(i-1)\epsilon\}< d(\alpha^- D^-,x_0) \leq i\epsilon\big\}\,,
\]
\[
\underline{\A}^+_j=
\big\{\alpha^+\in\Ga/\Ga_{D^+}:\max\{0,j\epsilon-1+\epsilon\} <
d(x_0,\alpha^+ D^+) \leq j\epsilon\big\}\,, \;\text{and}
\]
\[
\overline{\A}^+_j=\big\{\alpha^+ \in \Ga/\Ga_{D^+}:
\max\{0,j\epsilon-1\}< d(x_0,\alpha^+ D^+) \leq j\epsilon+2\epsilon\big\}\,.
\]

Let $\alpha^-\in\Ga/\Ga_{D^-}$ and $\alpha^+ \in\Ga/\Ga_{D^+}$ be such
that we have $d(\alpha^-D^-,\alpha^+D^+)>0$, $d(x_0, [\alpha^-D^-,
  \alpha^+D^+]) \leq \frac{\epsilon}{2}$ and $x_0\notin \alpha^-D^-
\cup\alpha^+D^+$.  Let $p_{\alpha^-,\alpha^+}$ be the closest point to
$x_0$ on the geodesic segment $[\alpha^-D^-, \alpha^+D^+]$ (see Figure
\ref{fig:perpcomm}).  Since $d(x_0, p_{\alpha^-, \alpha^+}) \leq
\frac{\epsilon}{2}$ and since closest point projections do not
increase the distances, we have
\[
d(\alpha^-D^-,x_0)-\frac{\epsilon}{2}\leq
d(\pi(u_{\alpha^-D^-,\,\alpha^+D^+}),x_0)-\frac{\epsilon}{2}\leq
d(\alpha^-D^-,p_{\alpha^-,\alpha^+}) \leq d(\alpha^-D^-,x_0)\,.
\]
Similarly, we have 
\[
d(x_0,\alpha^+D^+)-\frac{\epsilon}{2}\leq
d(p_{\alpha^-,\alpha^+},\alpha^+D^+) \leq d(x_0,\alpha^+D^+)\;.
\]  
Hence
\[
d(\alpha^-D^-,x_0)+d(x_0,\alpha^+D^+)-\epsilon\leq
d(\alpha^-D^-,\alpha^+D^+)
\leq d(\alpha^-D^-,x_0)+d(x_0,\alpha^+D^+)\,.
\]
Therefore
\begin{equation}\label{eq:mino}
  \text{if}\quad\alpha^- \in \A^-_i\quad\text{and}
  \quad\alpha^+\in\underline{\A}^+_{k-i}, \quad\text{then}
  \quad(\alpha^-,\alpha^+)\in\A_k\,.
\end{equation}
and
\begin{equation}\label{eq:majo}
  \text{if}\quad\alpha^- \in \A^-_{i}\quad\text{and}
  \quad(\alpha^-,\alpha^+)\in\A_{k}, \quad\text{then}
  \quad\alpha^+\in\overline{\A}^+_{k-i}\,.
\end{equation}

By integrating on the first factor the first formula in Theorem
\ref{theo:17ETDS2} applied by taking $A^+= \Ga\{x_0\}$, and by lifting
to $T^1\wt M$, for the weak-star convergence of measures on
$T^1_{x_0}\wt M$, we have (uniformly on $\epsilon$)
\begin{equation}\label{eq:prepminuside}
\lim_{s\,\epsilon\ra+\infty}\;\frac{\delta\;\|\mBM\|}
{e^{\delta\, s\,\epsilon}\;\|\sigma^+_{A^-}\|}
\sum_{\alpha^-\in\,\bigcup_{i=1}^s\A^-_i}\Delta_{v_{\alpha^-D^-,\{x_0\}}}
=\wt\sigma^-_{\{x_0\}}\,.
\end{equation}
For every $i\in \NN$, let us define
\[
a_i=\sum_{\alpha^-\in\A^-_i}\psi^-\big((u_{\alpha^-D^-,\{x_0\}})_-\big)\,.
\]
Since the negative endpoint map $p_- :T^1_{x_0}M\ra\partial_\infty\wt
M$ is a homeomorphism, its pushforward map $(p_-)_*$ on finite
measures is weak-star continuous. Note that we have the equality
$p_-(v_{\alpha^-D^-,\{x_0\}})= (u_{\alpha^-D^-,\{x_0\}})_-$. Hence by
evaluating on $\psi^-$ the image of Equation \eqref{eq:prepminuside}
by $(p_-)_*$, and by Equation \eqref{eq:skinningsingle}, we have
\begin{align}
\lim_{s\,\epsilon\ra+\infty}\;\;\frac{\delta\;\|\mBM\|}
    {e^{\delta\, s\,\epsilon}\;\|\sigma^+_{A^-}\|}
    \sum_{i=1}^s a_i &=\lim_{s\,\epsilon\ra+\infty}\;\;\frac{\delta\;\|\mBM\|}
    {e^{\delta\, s\,\epsilon}\;\|\sigma^+_{A^-}\|}
    \sum_{\alpha^-\,\in\,\bigcup_{i=1}^s\A^-_i}\psi^-((u_{\alpha^-D^-,\{x_0\}})_-)
    \nonumber\\&=(p_-)_*\wt\sigma^-_{\{x_0\}}(\psi^-)
    = \mu_{x_0}(\psi^-)\,.\label{eq:minuside}
\end{align}
Hence for every $\eta>0$, there exists $i_0=i_0(\eta,\epsilon)$ such
that if $s\geq i_0$, then
\begin{equation}\label{eq:minuside2}
  \frac{\|\sigma^+_{A^-}\|\;(\mu_{x_0}(\psi^-)-\eta)}
  {\delta\;\|\mBM\|}\;e^{\delta\, s\,\epsilon}\leq
  \sum_{i=1}^s a_i\leq \frac{\|\sigma^+_{A^-}\|\;(\mu_{x_0}(\psi^-)+\eta)}
  {\delta\;\|\mBM\|}\;e^{\delta\, s\,\epsilon}\,.
\end{equation}
%See also \cite[Equation (12)]{ParPau24ETDS}.
Note that $e^{2\,\delta\,\epsilon}-e^{-\delta}=(1-e^{-\delta})
e^{\bigO(\epsilon)}$. Similarly, now integrating on the second factor
the first formula in Theorem \ref{theo:17ETDS2} now applied with $A^-=
\Ga\{x_0\}$, taking the difference between the times $j\epsilon+
2\epsilon$ and $j\epsilon-1$, lifting to $T^1\wt M$, pushing forwards
by $(p_+)_*$ and evaluating on $\psi^+$, we have (uniformly on
$\epsilon$)
\begin{equation}\label{eq:pluside}
\lim_{j\,\epsilon\ra+\infty}\;\;\frac{\delta\;\|\mBM\|}
    {(1-e^{-\delta})\;e^{\delta\, j\,\epsilon}
      \;\|\sigma^-_{A^+}\|}\sum_{\alpha^+\in\overline{\A}^+_j}
    \psi^+\big((v_{\alpha^-D^-,\{x_0\}})_+\big) =
    e^{\bigO(\epsilon)}\,\mu_{x_0}(\psi^+)\,.
\end{equation}
With the above $\bigO(\epsilon)$, for every $\eta>0$, let $f_{\pm}
:[0,+\infty[\ra \RR$ be the smooth function
\[
f_{\pm}:s\mapsto \frac{(1-e^{-\delta})
\; e^{\bigO(\epsilon)}\;\|\sigma^-_{A^+}\|\;(\mu_{x_0}(\psi^+)\pm\eta)}
{\delta\;\|\mBM\|}\;e^{\delta\, (k-s)\,\epsilon}\,.
\]
It follows from Equation \eqref{eq:pluside} that there exists
$j_0=j_0(\eta,\epsilon)$ such that if $j\geq j_0$, then
\begin{equation}\label{eq:pluside2}
  f_{-}(k-j)\leq \sum_{\alpha^+\in\overline{\A}^+_j}
  \psi^+\big((v_{\alpha^-D^-,\{x_0\}})_+\big)\leq f_{+}(k-j)\,.
\end{equation}
Let us define
\[
C_1=\frac{\|\sigma^+_{A^-}\|\;(\mu_{x_0}(\psi^-)+\eta)}
  {\delta\;\|\mBM\|}
\quad\text{and}\quad C_2=\frac{(1-e^{-\delta}) \;e^{\bigO(\epsilon)} \;
  \|\sigma^-_{A^+}\| \;(\mu_{x_0}(\psi^+)+\eta)}{\delta\;\|\mBM\|}\,.
\]
By decomposing the sum defining $Z_k$ in Equation
\eqref{eq:defiSigmak} into thin slices of width $\epsilon$ of the
first index $\alpha^-$, and by using Equation \eqref{eq:majo}, we have
\begin{align}
  Z_k&=\sum_{i\in\NN}\;\sum_{\alpha^-\in\A^-_{i}}\Big(
  \psi^-\big((u_{\alpha^-D^-,\{x_0\}})_-\big)
\sum_{\substack{\alpha^+\in\Ga/\Ga_{D^+}\\(\alpha^-,\alpha^+)\in\A_{k}}}\;
\psi^+\big((v_{\{x_0\},\alpha^+D^+})_+\big)\Big)\int\psi^0\, ds
\nonumber\\&\leq \sum_{i\in\NN}\;a_i\Big(
\sum_{\alpha^+\in\overline{\A}^+_{k-i}}\;
\psi^+\big((v_{\{x_0\},\alpha^+D^+})_+\big)\Big)\int\psi^0\, ds\,.
\label{eq:controlZk}
\end{align}
Note that $\overline{\A}^+_{j}$ is empty if $j\leq -3$, hence
$\overline{\A}^+_{k-i}$ is empty if $i\geq k+3$.  We decompose the sum
$\sum_{i=0}^{k+2}$ into the sum $\sum_{i=0}^{k-j_0}$, where we can use
the upper bound in Equation \eqref{eq:pluside2} with $j=k-i$, and the
sum $\sum_{i=k-j_0+1}^{k+2}$, where we can use the estimate
$\sum_{i=k-j_0+1}^{k+3} a_i=\bigO(e^{\delta\, k\,\epsilon})$ following
from Equation \eqref{eq:minuside} and the estimate $\sum_{\alpha^+\in
  \overline{\A}^+_{k-i}} \psi^+\big( (v_{\alpha^-D^-,\{x_0\}})_+
\big)= \bigO(e^{\delta\, (k-i)\,\epsilon})$ following from Equation
\eqref{eq:pluside}. Hence, using Abel's summation formula and the fact
that $a_0=0$ for the equality below, we have
\begin{align*}
Z_k& \leq \sum_{i=0}^{k-j_0}\;a_i\;f_+(i)\int\psi^0\, ds+
\bigO_{j_0}(e^{\delta\, k\,\epsilon})\\&=
\Big(\Big(\sum_{i=0}^{k-j_0}\;a_i\Big)\;f_+(k-j_0)-\int_{u=0}^{k-j_0}
\Big(\sum_{i=0}^{u}\;a_i\Big)f'_+(u)\, du\Big)\int\psi^0\, ds+
\bigO_{j_0}(e^{\delta\, k\,\epsilon}) \,.
\end{align*}
Note that $f_+(k-j_0)=\bigO_{j_0}(1)$ and $f'_+(u)= -\,\delta\,\epsilon\;
C_2\;e^{\delta\, (k-u)\,\epsilon}=\bigO(e^{\delta\, k\,\epsilon})$ when
$u\leq k$. Using again $\sum_{i=1}^{k-j_0} a_i=\bigO(e^{\delta\,
  k\,\epsilon})$, subdividing the integral $\int_{u=0}^{k-j_0}$ into
$\int_{u=0}^{i_0}$ and $\int_{u=i_0}^{k-j_0}$ and using on the second
one the inequality $\sum_{i=1}^u a_i\leq C_1\;e^{\delta\,
  u\,\epsilon}$ following from the upper bound in Equation
\eqref{eq:minuside2}, we have
\begin{align*}
  Z_k&\leq \int_{u=i_0}^{k-j_0} (C_1\;e^{\delta u\epsilon})
  (\delta\,\epsilon\; C_2\;e^{\delta (k-u)\epsilon})\, du\int\psi^0\, ds+
  \bigO_{i_0,j_0}(e^{\delta k\epsilon})\nonumber\\&
  =\delta\,k\,\epsilon\,C_1\, C_2\;e^{\delta k\epsilon}\int\psi^0\, ds+
  \bigO_{i_0,j_0}(e^{\delta k\epsilon})\,.
\end{align*}
By the upper bound in Equation \eqref{eq:relatnu3tSigmak}, by
expliciting the values of $C_1$, $C_2$, by a geometric series
argument, and since $t=N\epsilon+\bigO(\epsilon)$, for $t$ large
enough, we hence have
\[
\frac{\delta\;\|\mBM\|}{t\;e^{\delta\,t}
  \;\|\sigma^+_{A^-}\|\;\|\sigma^{-}_{A^+}\|}\nu_{3,t}(\psi)\leq
\frac{e^{\bigO(\epsilon)}}{\|\mBM\|}
(\mu_{x_0}(\psi^-)+\eta)\;(\mu_{x_0}(\psi^+)+\eta)\;\int\psi^0\, ds+
\bigO_{i_0,j_0}\big(\frac{1}{t}\big)\;.
\]
By taking the upper limit as $t\ra+\infty$ and by letting $\eta\ra 0$,
we thus have
\[
\limsup_{t\ra+\infty} \frac{\delta\;\|\mBM\|}
{t\;e^{\delta\, t}\;\|\sigma^+_{A^-}\|\;\|\sigma^{-}_{A^+}\|}
\;\nu_{3,t}(\psi)\leq e^{\bigO(\epsilon)}
\frac{1}{\|\mBM\|}\,d\mu_{x_0}\otimes d\mu_{x_0}\otimes ds\;(\psi)\,.
\]

Since every element $v\in T^1\wt M$ of the support of $\psi$ satisfies
$d(\pi(v),x_0) \leq \frac{\epsilon}{2}$ and since by Equation
\eqref{eq:defiBuse} we have $|\,\beta_\xi(x,y)\,| \leq d(x,y)$ for all
$\xi\in \partial_\infty \wt M$ and $x,y\in\wt M$, by Equation
\eqref{eq:defimBM}, we have
\begin{equation}\label{eq:BMquasiprod}
d\mu_{x_0}\otimes d\mu_{x_0}\otimes ds\;(\psi)
=e^{\bigO(\epsilon)}\; \wtmBM(\psi)\,.
\end{equation}
Therefore by Equations \eqref{eq:nu1} and \eqref{eq:nu2}, we have
\[
\limsup_{t\ra+\infty} \frac{\delta\;\|\mBM\|}
{t\;e^{\delta\, t}\;\|\sigma^+_{A^-}\|\;\|\sigma^{-}_{A^+}\|}\;\nu_{1,t}(\psi)
\leq e^{\epsilon_3}\;\frac{\wtmBM(\psi)}{\|\mBM\|}\,.
\]
A similar argument replacing Equation \eqref{eq:majo} by Equation
\eqref{eq:mino}, and the upper bounds in Equations \eqref{eq:pluside2},
\eqref{eq:minuside2}, \eqref{eq:relatnu3tSigmak} by their lower bounds,
gives
\[
\liminf_{t\ra+\infty} \frac{\delta\;\|\mBM\|}
{t\;e^{\delta\, t}\;\|\sigma^+_{A^-}\|\;\|\sigma^{-}_{A^+}\|}\;\nu_{1,t}(\psi)
\geq e^{-\epsilon_3}\;\frac{\wtmBM(\psi)}{\|\mBM\|}\,.
\]
Let $\psi'$ be a continuous function with compact support on $T^1\wt
M$.  By covering its support with sets of the form $\pi^{-1}(B(x_0,
\frac{\epsilon}{2}))$ for finitely many $x_0\in T^1\wt M$, by using a
partition of unity, by a uniform approximation by product functions,
and by letting $\epsilon$ tends to $0$, we thus obtain that
\[
\lim_{t\ra+\infty} \frac{\delta\;\|\mBM\|}
{t\;e^{\delta\, t}\;\|\sigma^+_{A^-}\|\;\|\sigma^{-}_{A^+}\|}\;
\nu_{1,t}(\psi')= \frac{\wtmBM(\psi')}{\|\mBM\|}\,.
\]
This proves that the measures $\mu_t=\frac{\delta\;\|\mBM\|}
{t\;e^{\delta\, t}\;\|\sigma^+_{A^-}\|\;\|\sigma^{-}_{A^+}\|}
\sum_{\alpha\in\Perp(A^-,\,A^+,\,t)} \; \Leb_\alpha$ converge to
$\frac{\mBM}{\|\mBM\|}$ for the weak-star convergence of measures on
$T^1M$.

The total mass of the measure $\mu_t$ converges to $1$, since by
\cite[Theo.~1]{ParPau17ETDS} we have $\card\Perp(A^-,\,A^+,\,t)\sim
\frac{\|\sigma^+_{A^-}\|\;\|\sigma^{-}_{A^+}\|}{\delta\;\|\mBM\|}\;
e^{\delta\, t}$ as $t\ra+\infty$, and since by this exponential growth
property, most of the mass of the sum $\sum_{\alpha\in
  \Perp(A^-,\,A^+,\,t)} \; \Leb_\alpha$ is obtained when the length of
$\alpha$ is close to $t$.  Since $\frac{\mBM}{\|\mBM\|}$ is a
probability measure, there is hence no loss of mass in the above
weak-star convergence, and it is well-known that this implies the
narrow convergence of $\mu_t$ to $\frac{\mBM}{\|\mBM\|}$.  This
concludes the proof of the first claim of Theorem
\ref{theo:mainintro}.

\bigskip
The proof of the second claim of Theorem \ref{theo:mainintro} proceeds
very similarly, we only indicate the changes. Let $\ell$ and $\kappa$
be as in Theorem \ref{theo:17ETDS2}. Let $K$ be a compact subset of
$T^1M$, and $\wt K$ a compact subset of $T^1\wt M$ mapping to $K$. Let
$x_0\in \wt K$ and $\epsilon\in \interval[open]{0}{\frac{1}{2}}$.
Since $M$ is locally symmetric with finite volume, the geometric
compactification $\overline{\wt M}=\wt M\cup\partial_\infty\wt M$ is a
smooth manifold with boundary, and Hopf's parametrisation with respect
to $x_0$ is a smooth diffeomorphism. Furthermore, the
Patterson-Sullivan measure $\mu_{x_0}$ is up to a multiplicative
constant the smooth measure on $\partial_\infty\wt M$ invariant under
the group of isometries of $\wt M$ fixing $x_0$, and the
Bowen-Margulis measure on $T^1M$ is up to a multiplicative constant
the Liouville measure of $T^1M$.

We now start with $\psi\in C^\ell(T^1\wt M)$ with support in
$\pi^{-1}(B(x_0, \frac{\epsilon}{2}))$, which is, in Hopf's
parametrisation with respect to $x_0$, a product $(\xi,\eta,s) \mapsto
\psi^-(\xi)\;\psi^+(\eta)\;\psi^0(s)$ of three smooth functions
$\psi^-:\partial_\infty\wt M\ra\CC$, $\psi^+:\partial_\infty \wt
M\ra\CC$, $\psi^0:\RR\ra\CC$ with compact support. The Bowen-Margulis
measure is absolutely continuous with respect to the product measure
$d\mu_{x_0}\otimes d\mu_{x_0}\otimes ds$, with a smooth density, which
is bounded from above and from below by a positive constant on the
compact subset $\wt K$. Hence by Fubini's theorem, we have
\begin{equation}\label{eq:majoSobnorm}
\|\psi^-\|_{\ell}\;\|\psi^+\|_{\ell}\;\|\psi^0\|_{\ell}=\bigO_{\wt K}
(\|\psi\|_{\ell})\,.
\end{equation}
We extend $\psi^-$ and $\psi^+$ to smooth functions on $\wt M$ with
support in the complementary subset of $B(x_0,-\ln \epsilon)$, and
then on $T^1\wt M$ by being constant on the fibers of the footpoint
projection $\pi:T^1\wt M\ra\wt M$, with the same Sobolev norms up to a
multiplicative constant.

By the error term in Theorem \ref{theo:17ETDS2}, Equation
\eqref{eq:prepminuside} now becomes, for all $\phi\in
C^\ell(T^1_{x_0}\wt M)$ and $s>0$,
\[
\frac{\delta\;\|\mBM\|}
{e^{\delta\, s\,\epsilon}\;\|\sigma^+_{A^-}\|}
\sum_{\alpha^-\in\,\bigcup_{i=1}^s\A^-_i}\phi(v_{\alpha^-D^-,\{x_0\}})
=\wt\sigma^-_{\{x_0\}}(\phi)+\bigO(e^{-\kappa s \epsilon}\,\|\phi\|_{\ell})\,.
\]
Since $p_- :T^1_{x_0}M\ra\partial_\infty\wt M$ is now a smooth
diffeomorphism between compact manifolds, applying the above result to
$\phi=\psi^-\circ p_-$, with $C'_1=\frac {\|\sigma^+_{A^-}\|\;
  \mu_{x_0}(\psi^-)} {\delta\;\|\mBM\|}$, Equation \eqref{eq:minuside}
becomes
\begin{equation}\label{eq:minusideET}
  \sum_{i=1}^s a_i=C'_1\;e^{\delta\, s\,\epsilon}+
  \bigO(e^{(\delta-\kappa) s \epsilon}\,\|\psi^-\|_{\ell})\,.
\end{equation}

For appropriately choosen functions $\bigO(\cdot)$, and $C'_2=\frac
{(1-e^{-\delta})\;e^{\bigO(\epsilon)}\;\|\sigma^-_{A^+}\|\;
  \mu_{x_0}(\psi^+)}{\delta\;\|\mBM\|}$, we now define a smooth
function $f:\RR\ra \RR$ by
\[
f:s\mapsto C'_2\;e^{\delta(k-s)\epsilon}+
\bigO(e^{(\delta-\kappa)  (k-s) \epsilon}\,\|\psi^+\|_{\ell})\,,
\]
so that $f':s\mapsto-\,\delta\,\epsilon\,C'_2\;e^{\delta(k-s)\epsilon}
+ \bigO(\epsilon\,e^{(\delta-\kappa) (k-s) \epsilon}\,
\|\psi^+\|_{\ell})$.  Equations \eqref{eq:pluside} and
\eqref{eq:pluside2} become, for every $j\in\ZZ$,
\begin{equation}\label{eq:plusideET}
  \sum_{\alpha^+\in\overline{\A}^+_j}
    \psi^+\big((v_{\alpha^-D^-,\{x_0\}})_+\big)=f(k-j)\,.
  \end{equation}
Since $\overline{\A}^+_{k-i}$ is empty if $i\geq k+2$, Equation
\eqref{eq:controlZk} gives, by Abel's summation formula, that
\begin{align*}
  Z_k& \leq \sum_{i=0}^{k+2}\;a_i\;f(i)\int\psi^0\, ds
  \\&=\Big(\int\psi^0\, ds\Big)\;
  \left(\Big(\sum_{i=0}^{k+2}\;a_i\Big)\;f(k+2)-\int_{u=0}^{k+2}
  \Big(\sum_{i=0}^{u}\;a_i\Big)f'(u)\, du\right)\,.
\end{align*}
We have $f(k+2) = \bigO(\|\psi^+\|_{\ell})$ by the definition of $f$
and since by the Cauchy-Schwarz inequality, we have $\mu_{x_0}(\psi^+)
\leq \|\mu_{x_0}\|\|\psi^+\|_{0}\leq \|\mu_{x_0}\| \|\psi^+\|_{\ell}$.
Similarly, $\int\psi^0\, ds= \bigO_{\wt K}(\|\psi^0\|_{\ell})$ and by
Equation \eqref{eq:minusideET}, we have $\sum_{i=1}^{k+2} a_i=
\bigO(e^{\delta k\epsilon}\,\|\psi^-\|_{\ell})$. Hence, again by
Equation \eqref{eq:minusideET} and by the computation of $f'$,  we have
{\small
  \begin{align*}
  Z_k&\leq \bigO_{\wt K}(e^{\delta k\epsilon}\,\|\psi^-\|_{\ell}\;
  \|\psi^+\|_{\ell}\,\|\psi^0\|_{\ell})
\\&\quad+\int\psi^0\, ds
  \int_{u=0}^{k+2}\big(C'_1\;e^{\delta u\epsilon}+
  \bigO(e^{(\delta-\kappa) u \epsilon}\,\|\psi^-\|_{\ell})\big)
  \big(\delta\,\epsilon\,C'_2\;e^{\delta(k-u)\epsilon}+
  \bigO(\epsilon\,e^{(\delta-\kappa)  (k-u) \epsilon}\,\|\psi^+\|_{\ell})\big)\,du\,.
  \end{align*}
} $\!\!$By Equation \eqref {eq:majoSobnorm} and by expliciting the
values of $C'_1$, $C'_2$, we hence have
\begin{align*}
  Z_k&\leq \bigO_{\wt K}(e^{\delta k\epsilon}\,\|\psi\|_{\ell})+
  (k+2)\,\delta\,\epsilon\, C'_1\,C'_2\,e^{\delta k\epsilon} \int\psi^0\, ds
  \nonumber\\& \qquad+
  \bigO_{\wt K}\big(\int_{u=0}^{k+2}(e^{(\delta-\kappa)k\epsilon +\kappa u\epsilon}+
  e^{\delta k\epsilon -\kappa u\epsilon})
 \,\epsilon\,du\;\|\psi\|_{\ell}\big) \\
  &=\frac{(1-e^{-\delta})\,e^{\bigO(\epsilon)}\;
      \|\sigma^+_{A^-}\|\;\|\sigma^-_{A^+}\|}{\delta\;\|\mBM\|^2}
 \;k\,\epsilon \,e^{\delta k\epsilon }\;\mu_{x_0}(\psi^-)\,\mu_{x_0}(\psi^+)
 \int\psi^0\, ds+\bigO_{\wt K}(e^{\delta k\epsilon }\,\|\psi\|_{\ell})\;.
\end{align*}
Therefore, by the upper bound in Equation \eqref{eq:relatnu3tSigmak},
by a geometric series argument, and since
$t=N\epsilon+\bigO(\epsilon)$, we have
\begin{align*}
  \nu_{3,t}(\psi)&\leq \frac{e^{\bigO(\epsilon)}\;\|\sigma^+_{A^-}\|\;
    \|\sigma^-_{A^+}\|}{\delta\;\|\mBM\|}\;t\,e^{\delta t}\,
  \int \psi\;\frac{1}{\|\mBM\|} d\mu_{x_0}\otimes d\mu_{x_0} \otimes ds
+\bigO_{\wt K}(e^{\delta t}\,\|\psi\|_{\ell})\;.
\end{align*}
The same lower bound is proved by using the lower bound in Equation
\eqref{eq:relatnu3tSigmak} and by replacing Equation \eqref{eq:majo}
by Equation \eqref{eq:mino}. By Equations \eqref{eq:nu1},
\eqref{eq:nu2} and \eqref{eq:BMquasiprod}, we thus have
\begin{align*}
  \nu_{1,t}(\psi)&= \frac{e^{\epsilon_4}\;\|\sigma^+_{A^-}\|\;
    \|\sigma^-_{A^+}\|}{\delta\;\|\mBM\|}\;t\,e^{\delta t}\;
  \frac{\mBM(\psi)}{\|\mBM\|}\; +\bigO_{\wt K}(e^{\delta t}\,\|\psi\|_{\ell})\;.
\end{align*}
Let $\psi':T^1\wt M\ra\RR$ be a continuous function with support in
$\wt K$.  By covering $\wt K$ with sets of the form $\pi^{-1}(B(x_0,
\frac{\epsilon}{2}))$ for finitely many $x_0\in T^1\wt M$, by using a
smooth partition of unity and a smooth approximation by
product functions, and by letting $\epsilon$ tend to $0$, we thus
obtain the second claim of Theorem \ref{theo:mainintro}.

\medskip
This concludes the proof of Theorem \ref{theo:mainintro} in the
Introduction.

\section{Potentials and equidistribution of 
  common perpendiculars}
\label{sec:potential}

In this section, we prove Theorem \ref{theo:intropot}, after defining
the various objects that appear in its statement. We keep the notation
$\wt M$, $x_*$, $\Ga$, $M$, $\pi$, $(\flow{t})_{t\in\RR}$, $p_\pm$,
$v\mapsto v_\pm$ of Section \ref{sec:geometry}.

\subsection{Potentials and Gibbs measures}
\label{subsec:thermoform}

In this subsection, we briefly recall the thermodynamic formalism of
geodesic flows in negative curvature, referring to
\cite{PauPolSha15,BroParPau19, DilTho25} for precisions and further
developments.

Let $\wt F:T^1\wt M\ra \RR$ be a {\it potential} on $T^1\wt M$, that
is, a $\Ga$-invariant, bounded\footnote{See \cite[\S 3.2]{BroParPau19}
for a weakening of this assumption.} Hölder-continuous real-valued map
on $T^1\wt M$, and let $F:T^1M=\Ga\bs T^1\wt M\ra \RR$ be its induced
function.  For all $x,y\in \wt M$, let us define the {\it amplitude}
of $\wt F$ between $x$ and $y$ to be $\int_x^y\wt F=0$ if $x=y$ and
otherwise $\int_x^y\wt F= \int_{0}^{d(x,y)} \wt F(\flow t v) \;dt$
where $v$ is the tangent vector at $x$ to the geodesic segment from
$x$ to $y$. If $\alpha:[a,b]\ra M$ is a locally geodesic segment in
the orbifold $M$, and if $\wt \alpha:[a,b]\ra \wt M$ is any lift of
$\alpha$, we define the {\it amplitude} of $F$ along $\alpha$ to be
$\int_\alpha F=\int_{\wt \alpha}\wt F= \int_{\wt \alpha(a)}^{\wt
  \alpha(b)}\wt F$.

The {\it critical exponent} of $F$ is the weighted (by the exponential
amplitudes of $F$) orbital growth rate of the group $\Ga$, defined by
\[
\delta_F= \lim_{n\ra+\infty}\;\frac{1}{n}\;\ln\;\Big(
    \sum_{\ga\in\Ga,\;n-1< d(x_*,\ga x_*)\leq n} \;\;
    \exp\big(\int_{x_*}^{\ga x_*} \wt F\;\big)\Big)\;.
\]
This limit exists and is independent of the choice of $x_*$. We have
$\delta_F\in \;\interval[open]{-\infty}{+\infty}$, $\delta_{F+c}=
\delta_{F}+c$ for every constant $c\in\RR$ and $\delta_{F\circ\iota}=
\delta_F$ for $\iota: v\mapsto -\,v$ the opposite map on $T^1\wt M$.
See \cite[Chap.~4, 6]{PauPolSha15} for equivalent definitions (in
particular for the equality with the topological pressure of $F$).

The (normalised) {\it Gibbs cocycle} of the potential $\wt F$ (as for
instance defined by Hamenstädt) is the function $C^F:\partial_\infty
\wt M\times \wt M\times \wt M\ra \RR$, defined by the following limit
of differences of amplitudes for the renormalised potential
\[
(\xi,x,y)\mapsto C^F_\xi(x,y)= \lim_{t\ra+\infty}
\int_y^{\xi_t}(\wt F-\delta_F)-\int_{x}^{\xi_t}(\wt F -\delta_F)\;,
\]
where $t\mapsto \xi_t$ is any geodesic ray in $\wt M$ converging to
$\xi$.  The Gibbs cocycle is $\Ga$-invariant (for the diagonal action)
and locally Hölder-continuous.

A {\it Patterson density} for $(\Ga,\wt F)$ is a $\Ga$-equiv\-ariant
family $(\mu^F_{x})_{x\in \wt M}$ of pairwise absolutely continuous
(positive, Borel) measures on $\partial_\infty \wt M$, whose support
is $\Lambda\Ga$, such that
\[
\ga_*\mu^F_x=\mu^F_{\ga x}{\rm ~~~and~~~}
\frac{d\mu^F_x}{d\mu^F_y}(\xi) = e^{-C^F_\xi(x,\,y)}
\]
for all $\ga\in\Ga$ and $x,y\in \wt M$, and for (almost) every
$\xi\in\partial_\infty \wt M$. Patterson's classical construction
gives its existence, see \cite[Prop.~3.9]{PauPolSha15}.

The {\it Gibbs measure} (or {\it equilibrium state}) on $T^1\wt M$
associated with a pair $(\mu^{F\circ\iota}_{x})_{x\in \wt M}$ and
$(\mu^F_{x})_{x\in \wt M}$ of Patterson densities for $(\Ga,\wt F\circ
\iota)$ and $(\Ga,\wt F)$ is the $\sigma$-finite nonzero measure $\wt
m_F$ on $T^1\wt M$ defined in \cite[Eq.~(43)]{PauPolSha15}, using the
Hopf parametrisation $v\mapsto (v_-,v_+,s)$ with respect to $x_*$, by
\begin{equation}\label{eq:defiGibbs}
d\wt m_F(v)=
e^{C^{F\circ\iota}_{v_-}(x_*,\,\pi(v))+C^F_{v_+}(x_*,\,\pi(v))}
\;d\mu^{F\circ\iota}_{x_*}(v_-)\;d\mu^F_{x_*}(v_+)\;ds\;.
\end{equation}
It is independent of the choice of $x_*$. It is $\Ga$-invariant and
$(\flow t)_{t\in\RR}$-invariant.  Therefore it induces\footnote{See
for instance \cite[\S 2.6]{PauPolSha15} for details on the definition
of the induced measure since $\Ga$ might have torsion, hence it does
not act freely on $T^1\wt M$.} a $\sigma$-finite nonzero $(\flow
t)_{t\in\RR}$-invariant measure on $T^1M=\Ga\bs T^1\wt M$, called the
{\it Gibbs measure} on $\Ga\bs T^1\wt M$ for the potential $F$ and
denoted by $m_F$. If $m_F$ is finite, then the above Patterson
densities and the Gibbs measure are unique up to a scalar multiple,
see \cite[\S 5.3]{PauPolSha15}.

Let $D$ be a nonempty proper closed convex subset of $\wt M$, and let
$A=\Ga D$ be its image in $M=\Ga\bs \wt M$. Using the endpoint
homeomorphisms $p_\pm:v\mapsto v_\pm$ from $\normalpm {D}$ to
$\partial_{\infty}\wt M -\partial_{\infty}D$, the inner and outer {\it
  skinning measures} $\wt \sigma^-_{D}$ on $\normalin{D}$ and $\wt
\sigma^+_{D}$ on $\normalout{D}$ associated with the Patterson
densities $(\mu^{F\circ\iota}_{x})_{x\in \wt M}$ and
$(\mu^F_{x})_{x\in \wt M}$ for $(\Ga,\wt F\circ \iota)$ and $(\Ga,\wt
F)$ respectively, are defined in \cite[Eq.~(7.1)]{BroParPau19} by
\begin{equation}\label{eq:defskin}
d\,\wt\sigma^-_D(v) =
e^{C^{F\circ\iota}_{v_-}(x_*,\,\pi(v))}\;d\mu^{F\circ\iota}_{x_*}(v_-)
\quad\text{and}\quad d\,\wt\sigma^+_D(v) =
e^{C^F_{v_+}(x_*,\,\pi(v))}\;d\mu^F_{x_*}(v_+) \;.
\end{equation}
They do not depend on $x_*$. If the $\Ga$-equivariant family $\D=(\ga
D)_{\ga\in\Ga/\Ga_D}$ is locally finite, the {\it inner and outer
  skinning measures} $\wt\sigma^-_{\D}$ and $\wt\sigma^+_{\D}$ of $D$
on $T^1\wt M$ are the $\Ga$-invariant locally finite measures on
$T^1\wt M$ defined by
\[
\wt\sigma^\pm_{\D}=\sum_{\ga\in\Ga/\Ga_D}\wt\sigma^\pm_{\ga D}\;.
\]
They induce locally finite measures on $T^1 M$, denoted by
$\sigma^-_{A}$ and $\sigma^+_{A}$, called the {\it inner and outer
  skinning measures} of $A$ on $T^1 M$ associated with the potential
$F$.

\subsection{Proof of Theorem \ref{theo:intropot}}
\label{subsec:proof2}

The proof proceeds similarly to the proof of the first claim of
Theorem \ref{theo:mainintro}, we only indicate the modifications. Let
$D^\pm,A^\pm$ be as in the beginning of Section \ref{sec:proof}. We
will use the following analog of Theorem \ref{theo:17ETDS2}, see
\cite[Theo.~1.4]{BroParPau19}.

\btheo\label{theo:19book1_4} Assume that $\delta_F>0$, that $m_F$ is
finite and mixing for the geodesic flow of $M$, and that
$\sigma^+_{A^-}$ and $\sigma^{-}_{A^+}$ are finite and nonzero. For
the weak-star convergence of measures on $T^1M\times T^1M$, we have
\[
\;\;\;\lim_{t\ra+\infty}\; \delta_F\;\|m_F\|\;e^{-\delta_F\, t}
\sum_{\alpha\in\Perp(A^-,\,A^+,\,t)} \; e^{\int_\alpha F}\;\Delta_{v^-_\alpha}
\otimes\Delta_{v^+_\alpha}\;=\; \sigma^+_{A^-}\otimes
\sigma^-_{A^+}\,. \quad\Box
\]
\etheo

Let $\epsilon,x_0,\psi,\psi^\pm,\psi^0$ be as in the proof of the
first claim of Theorem \ref{theo:mainintro}. Let $t>1$ and
$N=\lceil\frac{t}{\epsilon}\rceil$.  We now consider the
$\Ga$-invariant Borel measures on $T^1\wt M$ defined by
\[
\nu_{1,t}=\sum_{\stackrel{\alpha^-\in\Ga/\Ga_{D^-},\,
  \alpha^+\in\Ga/\Ga_{D^+}}{0<d( \alpha^-D^-,\,\alpha^+D^+)\leq t}}
e^{\int_{[\alpha^-D^-,\,\alpha^+D^+]}\wt F}\;\Leb_{[\alpha^-D^-,\alpha^+D^+]}
\quad\text{and}
\]
\[
\nu_{2,t}=\sum_{\stackrel{\alpha^-\in\Ga/\Ga_{D^-},\,
  \alpha^+\in\Ga/\Ga_{D^+}}{0<d( \alpha^-D^-,\,\alpha^+D^+)\leq t}}
e^{\int_{[\alpha^-D^-,\,\alpha^+D^+]}\wt F}\;\Delta_{u_{\alpha^-D^-,\alpha^+D^+}}\otimes
\Delta_{v_{\alpha^-D^-,\alpha^+D^+}}\otimes ds\,.
\]
Let $\alpha^-\in\Ga/\Ga_{D^-}$ and $\alpha^+\in\Ga/\Ga_{D^+}$ such
that $x_0 \notin \alpha^-D^- \cup\alpha^+D^+$ and $d(\alpha^-D^-,
\alpha^+D^+)>0$. We keep the notation of Figure \ref{fig:perpcomm} in
the proof of Theorem \ref{theo:mainintro}. By the additivity of the
amplitudes, we have
\[
\int_{[\alpha^-D^-,\,\alpha^+D^+]}\wt F=
  \int_{\pi(u_{\alpha^-D^-,\,\alpha^+D^+})}^{p_{\alpha^-,\,\alpha^+}}\wt F+
\int_{p_{\alpha^-,\,\alpha^+}}^{\pi(v_{\alpha^-D^-,\,\alpha^+D^+})}\wt F\,.
\]
Recall that $d(x_0,p_{\alpha^-,\,\alpha^+})\leq\frac{\epsilon}{2}$ and
that $d(\pi(u_{\alpha^-D^-,\,\alpha^+D^+}), \pi(u_{\alpha^-D^-,
  \,\{x_0\}})) \leq \frac{\epsilon}{2}$ since closest point maps do
not increase the distance. Hence by \cite[Lem.~3.2]{PauPolSha15},
since $\wt F$ is bounded, there exists a constant $c_2\in \;]0,1]$
such that if $\epsilon$ is small enough, we have
\[
\Big|\;\int_{\pi(u_{\alpha^-D^-,\,\{x_0\}})}^{x_0} \wt F-
\int_{\pi(u_{\alpha^-D^-,\,\alpha^+D^+})}^{p_{\alpha^-,\,\alpha^+}}\wt F \;\Big|
=\bigO(\epsilon^{c_2})\,,
\]
and similarly
\[
\Big|\;\int_{x_0}^{\pi(v_{\{x_0\},\alpha^+D^+})} \wt F-
\int_{p_{\alpha^-,\,\alpha^+}}^{\pi(v_{\alpha^-D^-,\,\alpha^+D^+})}\wt F \;\Big|
=\bigO(\epsilon^{c_2})\,.
\]
Hence with $\nu_{3,t}$ the $\Ga$-invariant Borel measure on $T^\wt M$
now defined by
\[
\nu_{3,t}=\sum_{\stackrel{\alpha^-\in\Ga/\Ga_{D^-},\,
    \alpha^+\in\Ga/\Ga_{D^+}}{\stackrel{x_0\,\notin\, \alpha^-D^-
\cup\alpha^+D^+}{0<d( \alpha^-D^-,\,\alpha^+D^+)\leq t}}}
e^{\int_{[\alpha^-D^-,\,\{x_0\}]}\wt F\;+\;\int_{[\{x_0\},\,\alpha^+D^+]}\wt F}\;
\Delta_{(u_{\alpha^-D^-,\{x_0\}})_-}\otimes
\Delta_{(v_{\{x_0\},\alpha^+D^+})_+}\otimes ds\,,
\]
Equations \eqref{eq:nu1} and \eqref{eq:nu2} are still satisfied.

For all $i,j,k\in\NN$, we define $\A_{k},\A^-_i,\underline{\A}^+_j,
\overline{\A}^+_j$ exactly as previously, so that Equations
\eqref{eq:mino} and \eqref{eq:majo} are still satisfied. Now, let
\begin{align*}
Z_k=\sum_{(\alpha^-,\alpha^+)\in\A_{k}} \;&e^{\int_{[\alpha^-D^-,\,\{x_0\}]}\wt F}\;
\psi^-\big((u_{\alpha^-D^-,\{x_0\}})_-\big) \times
\\&e^{\int_{[\{x_0\},\,\alpha^+D^+]}\wt F}
\;\psi^+\big((v_{\{x_0\},\alpha^+D^+})_+\big)\int\psi^0\, ds\,,
\end{align*}
so that Equation \eqref{eq:relatnu3tSigmak} is still satisfied. Using
Theorem \ref{theo:19book1_4} with $A^+= \Ga\{x_0\}$ instead of Theorem
\ref{theo:17ETDS2}, Equation \eqref{eq:prepminuside} becomes
\[
\lim_{s\,\epsilon\ra+\infty}\;\frac{\delta_F\;\|m_F\|}
{e^{\delta_F s\,\epsilon}\;\|\sigma^+_{A^-}\|}
\sum_{\alpha^-\in\,\bigcup_{i=1}^s\A^-_i}e^{\int_{[\alpha^-D^-,\,\{x_0\}]}\wt F}\;
\Delta_{v_{\alpha^-D^-,\{x_0\}}} =\wt\sigma^-_{\{x_0\}}\,.
\]
For every $i\in \NN$, let us now define
\[
a_i=\sum_{\alpha^-\in\A^-_i}e^{\int_{[\alpha^-D^-,\,\{x_0\}]}\wt F}\;
\psi^-\big((u_{\alpha^-D^-,\{x_0\}})_-\big)\,.
\]
By the left hand side of Equation \eqref{eq:defskin}, we have
$(p_-)_*\wt\sigma^-_{\{x_0\}}= \mu_{x_0}^{F\circ\iota}$, hence
Equation \eqref{eq:minuside} becomes
\[
\lim_{s\,\epsilon\ra+\infty}\;\;\frac{\delta_F\;\|m_F\|}
    {e^{\delta_F s\,\epsilon}\;\|\sigma^+_{A^-}\|}
    \sum_{i=1}^s a_i = \mu^{F\circ\iota}_{x_0}(\psi^-)\,.
\]
Similarly, using Theorem \ref{theo:19book1_4} with $A^-= \Ga\{x_0\}$
instead of Theorem \ref{theo:17ETDS2}, since $(p_+)_*\wt
\sigma^+_{\{x_0\}}= \mu_{x_0}^{F}$ by the right hand side of Equation
\eqref{eq:defskin}, Equation \eqref{eq:pluside} becomes
\[
\lim_{j\,\epsilon\ra+\infty}\;\;\frac{\delta_F\;\|m_F\|}
    {(1-e^{-\delta_F})\;e^{\delta_F j\,\epsilon}
      \;\|\sigma^-_{A^+}\|}\sum_{\alpha^+\in\overline{\A}^+_j}
    e^{\int_{[\{x_0\},\,\alpha^+D^+]}\wt F}\psi^+\big((v_{\alpha^-D^-,\{x_0\}})_+\big)
    = e^{\bigO(\epsilon)}\,\mu^F_{x_0}(\psi^+)\,.
\]
We now define
\[
f_{\pm}:s\mapsto \frac{(1-e^{-\delta_F})
\; e^{\bigO(\epsilon)}\;\|\sigma^-_{A^+}\|\;(\mu^F_{x_0}(\psi^+)\pm\eta)}
{\delta_F\;\|m_F\|}\;e^{\delta_F (k-s)\epsilon}\,,
\]
\[
C_1=\frac{\|\sigma^+_{A^-}\|\;(\mu^{F\circ\iota}_{x_0}(\psi^-)+\eta)}
  {\delta_F\;\|m_F\|}
\quad\text{and}\quad C_2=\frac{(1-e^{-\delta_F}) \;e^{\bigO(\epsilon)} \;
  \|\sigma^-_{A^+}\| \;(\mu^F_{x_0}(\psi^+)+\eta)}{\delta_F\;\|m_F\|}\,.
\]
As in the proof of the first claim of Theorem \ref{theo:mainintro},
we have
\begin{align*}
  Z_k&\leq \sum_{i\in\NN}\;a_i\Big(
\sum_{\alpha^+\in\overline{\A}^+_{k-i}}\;e^{\int_{[\{x_0\},\,\alpha^+D^+]}\wt F}
\psi^+\big((v_{\{x_0\},\alpha^+D^+})_+\big)\Big)\int\psi^0\, ds
\\ &\leq\delta_F\,k\,\epsilon\,C_1\, C_2\;e^{\delta_F k\epsilon}\int\psi^0\, ds+
  \bigO_{i_0,j_0}(e^{\delta_F k\epsilon})\,,
\end{align*}
\[
\text{and}\quad
\limsup_{t\ra+\infty} \frac{\delta_F\;\|m_F\|} {t\;e^{\delta_F
    t}\;\|\sigma^+_{A^-}\|\;\|\sigma^{-}_{A^+}\|}
\;\nu_{3,t}(\psi)\leq e^{\bigO(\epsilon)}
\frac{1}{\|m_F\|}\,d\mu^{F\circ\iota}_{x_0}\otimes
d\mu^F_{x_0}\otimes ds\;(\psi)\,.\quad
\]
By \cite[Lem.~3.4 (i)]{PauPolSha15}, since $\wt F$ is bounded, there
exists a constant $c_2\in \;]0,1]$ such that for all
$\xi\in\partial_\infty\wt M$ and $x,y\in\wt M$ with $d(x,y)\leq 1$,
we have
\[
\max\big\{\big|\,C^{F\circ \iota}_\xi(x,y)\,\big|,\;\big|\,C^F_\xi(x,y)\,
\big|\big\}=\bigO(d(x,y)^{c_2})\,.
\]
Since every element $v\in T^1\wt M$ of the support of $\psi$ satisfies
$d(\pi(v),x_0) \leq \frac{\epsilon}{2}$ and by Equation
\eqref{eq:defiGibbs} with $x_*=x_0$, we thus have
\begin{equation}\label{eq:Gibbsquasiprod}
d\mu^{F\circ\iota}_{x_0}\otimes d\mu^F_{x_0}\otimes ds\;(\psi)
=e^{\bigO(\epsilon^{c_2})}\; \wt m_F(\psi)\,.
\end{equation}
The remainder of the proof of Theorem \ref{theo:19book1_4} proceeds as
the one of the first claim of Theorem \ref{theo:mainintro}, using
\cite[Theo.~1.5 (1)]{BroParPau19} instead of
\cite[Theo.~1]{ParPau17ETDS} to pass from the weak-star convergence
to the narrow convergence. This concludes the proof of Theorem
\ref{theo:intropot}.

{\small \bibliography{../biblio}

\begin{thebibliography}{BAPP}

\bibitem[Bab]{Babillot02}
M.~Babillot.
\newblock {\it On the mixing property for hyperbolic systems}.
\newblock {Israel J. Math. {\bf 129} (2002), 61--76}.

\bibitem[Bourb]{Bourbaki71a}
N.~Bourbaki.
\newblock {\it Topologie g\'en\'erale}.
\newblock {chap.~1 \`a 4, Hermann, 1971}.

\bibitem[Bourd]{Bourdon95}
M.~Bourdon.
\newblock {\it Structure conforme au bord et flot g\'eod\'esique
  d'un CAT$(-1)$ espace}.
\newblock {L'Ens. Math. {\bf 41} (1995) 63--102}.

\bibitem[BrH]{BriHae99}
M.~R.~Bridson and A.~Haefliger.
\newblock {\it Metric spaces of non-positive curva\-tu\-re}.
\newblock {Grund. math. Wiss. {\bf 319}, Springer Verlag, 1999}.

\bibitem[BrPP]{BroParPau19}
A.~Broise-Alamichel, J.~Parkkonen, and F.~Paulin.
\newblock {\it Equidistribution and counting under equilibrium states
  in negative curvature and trees. Applications to non-Archimedean
  Diophantine approximation}.
\newblock {With an Appendix by J.~Buzzi. Prog. Math. {\bf 329},
  Birkhäuser, 2019}.

\bibitem[DT]{DilTho25}
  C. Dilsavor and D.~J.~Thompson.
\newblock {\it Gibbs measures for geodesic flow on CAT(-1) spaces}.
\newblock {Preprint {\tt [arXiv:2309.03297]}}.
  
\bibitem[EPP]{ErlParPau24}
V.~Erlandsson, J.~Parkkonen, and F.~Paulin.
\newblock {\it Counting and equidistribution of reciprocal
  closed geodesic in negative curvature}.
\newblock {In preparation}.

\bibitem[LP]{LiPan22} 
J.~Li and W.~Pan.
\newblock {\it Exponential mixing of geodesic flows for geometrically 
finite hyperbolic manifolds with cusps}.
\newblock {Invent. Math. {\bf 231} (2022) 931--1021}.

\bibitem[LS]{LiSta23} 
X.~Li and B.~Staffa.
\newblock {\it On the equidistribution of closed geodesics and geodesic nets}.
\newblock Trans. Amer. Math. Soc. {\bf 376} (2023) 8825--8855.

\bibitem[OhS1]{OhSha12}
H.~Oh and N.~Shah.
\newblock {\it The asymptotic distribution of circles in the
orbits of Kleinian  groups}.
\newblock {Invent. Math. {\bf 187} (2012) 1--35}.

\bibitem[OhS2]{OhSha13}
H.~Oh and N.~Shah.
\newblock {\it Equidistribution and counting for orbits of 
geometrically finite hyperbolic groups}.
\newblock {J. Amer. Math. Soc. {\bf 26} (2013) 511--562}.

\bibitem[OP]{OtaPei04}
J.-P.~Otal and M.~Peign\'e.
\newblock {\it Principe variationnel et groupes kleiniens}.
\newblock {Duke Math. J. {\bf 125} (2004) 15--44}.

\bibitem[PaP1]{ParPau14ETDS}
J.~Parkkonen and F.~Paulin.
\newblock {\it Skinning measure in negative curvature and
  equidistribution of equidistant submanifolds}.
\newblock {Erg. Theo. Dyn. Sys. {\bf 34} (2014) 1310--1342}.

\bibitem[PaP2]{ParPau16LMS}
J.~Parkkonen and F.~Paulin.
\newblock {\it Counting arcs in negative curvature}.
\newblock {In "Geometry, Topology and Dynamics in Negative Curvature" 
(ICM 2010 satellite conference, Bangalore), C.~S.~Aravinda, 
T.~Farrell, J.-F.~Lafont eds, London Math. Soc. Lect. Notes 
{\bf 425}, Cambridge Univ. Press, 2016}.

\bibitem[PaP3]{ParPau17ETDS}
J.~Parkkonen and F.~Paulin.
\newblock {\it Counting common perpendicular arcs in negative
  curvature}.
\newblock {Erg. Theo. Dyn. Sys. {\bf 37} (2017) 900--938}.

\bibitem[PaP4]{ParPau17MA}
J.~Parkkonen and F.~Paulin.
\newblock {\it Counting and equidistribution in Heisenberg groups}.
\newblock {Math. Annalen {\bf 367} (2017) 81--119}.

\bibitem[PaP5]{ParPau17CIRM}
J.~Parkkonen and F.~Paulin.
\newblock {\it A survey of some arithmetic applications of 
ergodic theory in negative curvature}.
\newblock {In "Ergodic theory and negative curvature" CIRM Jean
Morley Chair subseries, B.~Hasselblatt ed, Notes Math. 2164,
pp. 293--326, Springer Verlag, 2017}.

\bibitem[PaP6]{ParPau22MPCPS}
J.~Parkkonen and F.~Paulin.
\newblock {\it Counting and equidistribution 
in quaternionic Heisenberg groups}.
\newblock {Math. Proc. Cambridge Phil. Soc. {\bf 173} (2022) 67--104}.

\bibitem[PaP7]{ParPau24ETDS}
J.~Parkkonen and F.~Paulin.
\newblock {\it Joint partial equidistribution of Farey rays
  in negatively curved manifolds and trees}.
\newblock{Erg. Theo. Dyn. Sys {\bf 44} (2024) 2700--2736}.

\bibitem[PaP8]{ParPau25}
J.~Parkkonen and F.~Paulin.
\newblock {\it Divergent geodesics, ambiguous closed geodesics and
  the binary additive divisor problem}.
\newblock \newblock {Preprint {\tt [arXiv:2409.18251]}}.

\bibitem[ParPS]{ParPauSay25}
J.~Parkkonen, F.~Paulin and R.~Sayous.
\newblock {\it Equidistribution of divergent geodesic in
  negative curvature}.
\newblock {In preparation}.

\bibitem[PauPS]{PauPolSha15}
F.~Paulin, M.~Pollicott, and B.~Schapira.
\newblock {\it Equilibrium states in negative curvature}.
\newblock {Astérisque {\bf 373}, Soc. Math. France, 2015}.

\bibitem[PolWar]{PolWar}
M.~Pollicott and K.~War.
\newblock {\it Counting geodesic loops on surfaces of genus at least 2 without conjugate points}.
\newblock \newblock {Preprint {\tt [arXiv:2309.14099]}}.

\bibitem[Rob]{Roblin03}
T.~Roblin.
\newblock {\it Ergodicit\'e et \'equidistribution en courbure n\'egative}.
\newblock {M\'emoire Soc. Math. France, {\bf 95} (2003)}.


\end{thebibliography}
%\bibliography{viitteet} 
}

\bigskip
{\small\noindent \begin{tabular}{l} 
Department of Mathematics and Statistics, P.O. Box 35\\ 
40014 University of Jyv\"askyl\"a, FINLAND.\\
{\it jouni.t.parkkonen@jyu.fi}
\end{tabular}

\medskip\noindent \begin{tabular}{l}
Laboratoire de mathématique d'Orsay, UMR 8628 CNRS,\\
Universit\'e Paris-Saclay, 91405 ORSAY Cedex, FRANCE\\
{\it frederic.paulin@universite-paris-saclay.fr}
\end{tabular}}

\end{document}